\documentclass[10pt]{amsart}
\usepackage[latin1]{inputenc}

\usepackage[mathcal]{eucal}
\usepackage{amsmath}
\usepackage{amsfonts}
\usepackage{amssymb}
\usepackage{amsthm}
\usepackage{graphicx,epsfig}
\usepackage{amscd}
\DeclareMathOperator{\re}{Re}
\DeclareMathOperator{\im}{Im}

\newcommand{\nil}{\mathrm{Nil}_3}
\newcommand{\vp}{\varphi}
\newcommand{\D}{\mathbb{D}}

\newcommand{\R}{\mathbb{R}}
\newcommand{\h}{\mathbb{H}^2}
\newcommand{\Z}{\mathbb{Z}}
\newcommand{\N}{\mathbb{N}}
\newcommand{\C}{\mathbb{C}}
\newcommand{\s}{\mathbb{S}}
\newcommand{\LL}{\mathbb{L}}

\newcommand{\om}{\omega}
\newcommand{\si}{\sigma}

\newcommand{\rmT}{\mathrm{T}}

\newcommand{\rmd}{\mathrm{d}}
\newcommand{\rmO}{\mathrm{O}}

\newcommand{\cH}{{\mathcal H}}

\newcommand{\cC}{{\mathcal C}}
\newcommand{\cD}{{\mathcal D}}

\begin{document}

\newtheorem{thm}{Theorem}[section]
\newtheorem*{thmintro}{Theorem}
\newtheorem{cor}[thm]{Corollary}
\newtheorem{prop}[thm]{Proposition}
\newtheorem{app}[thm]{Application}
\newtheorem{lemma}[thm]{Lemma}
\newtheorem*{lemmaintro}{Lemma}
\newtheorem{notation}[thm]{Notations}
\newtheorem{hypothesis}[thm]{Hypothesis}

\newtheorem{defin}[thm]{Definition}
\newenvironment{defn}{\begin{defin} \rm}{\end{defin}}
\newtheorem{remk}[thm]{Remark}
\newenvironment{rem}{\begin{remk} \rm}{\end{remk}}
\newtheorem{exa}[thm]{Example}
\newenvironment{ex}{\begin{exa} \rm}{\end{exa}}

\title[Half-space theorem]{Half-space theorem, embedded minimal annuli and minimal graphs in the Heisenberg group}
\author{Beno\^\i t Daniel and Laurent Hauswirth}
\date{}


\address{Universit\'e Paris 12, D\'epartement de Math\'ematiques, UFR des Sciences et Technologies, 61 
avenue du G\'en\'eral de Gaulle, B\^at. P3, 4e \'etage, 94010 Cr\'eteil cedex, FRANCE}
\email{daniel@univ-paris12.fr}
\address{Universit\'e de Marne-la-Vall\'ee, D\'epartement de Math\'ematiques, Cit\'e Descartes, 5 bd Descartes, Champs-sur-Marne, 77454 Marne-la-Vall\'ee, Cedex 2, FRANCE}
\email{hauswirth@univ-mlv.fr}

\begin{abstract}
We construct a one-parameter family of properly embedded minimal annuli in the Heisenberg group $\nil$ endowed 
with a left-invariant Riemannian metric. These annuli are not rotationally invariant. This family gives a 
vertical half-space theorem and proves that each complete minimal
graph in $\nil$ is entire. Also, the sister surface of an entire minimal
graph in $\nil$ is an entire constant mean curvature $\frac12$ graph in $\h \times \R$, and conversely.
This gives a classification of all entire constant mean curvature $\frac12$ graphs in $\h \times \R$.
Finally we construct properly embedded constant mean curvature $\frac12$ annuli in $\h \times \R$.
\end{abstract}

\maketitle

\section{Introduction}

This paper deals with
global properties of minimal and constant mean curvature (CMC) surfaces in Riemannian homogeneous manifolds. 
Some interesting properties are the existence of a Hopf-type holomorphic quadratic differential (see \cite{AR1}, \cite{AR2}) and Lawson-type local isometric correspondences, in particular between minimal surfaces in the Heisenberg 
group $\nil$ endowed with a left-invariant Riemannian metric 
and CMC $\frac12$ surfaces in $\h \times \R$ (see \cite{dan}). Two surfaces related by this correspondence are 
called sister surfaces.


In this paper we first construct a one-parameter family of properly embedded minimal annuli that are somewhat
analogous to catenoids of $\R^3$ with a ``horizontal axis".

\begin{thmintro}[see theorem \ref{propertieshorcat}]
There exists a one-parameter family $(\cC_{\alpha})_{\alpha>0}$ of properly embedded minimal annuli in $\nil$, called ``horizontal catenoids", having the following properties:
\begin{itemize}
\item the annulus $\cC_{\alpha}$ is \emph{not} invariant by a one-parameter group of isometries,
\item the intesection of $\cC_{\alpha}$ and any vertical plane of equation $x_2=c$ ($c\in\R$) is a non-empty closed embedded convex curve,
\item the annulus $\cC_{\alpha}$ is invariant by rotations of angle $\pi$ around the $x_1$, $x_2$ and $x_3$ axes and the $x_2$-axis is contained in the ``interior'' of $\cC_{\alpha}$,
\item the annulus $\cC_{\alpha}$ is conformally equivalent to $\C\setminus\{0\}$.
\end{itemize}
(The model we use for $\nil$ is described in section \ref{preliminaries}.)
\end{thmintro}

Up to now, the only known examples of complete minimal surfaces in $\nil$ were surfaces invariant by a one-parameter group of isometries \cite{mercuri}, periodic surfaces \cite{AR2} and entire graphs \cite{dan2}. The annuli we construct are the first non-trivial examples of annuli; they are very different from the rotationally invariant catenoids, which are of hyperbolic conformal type. The existence of annuli of parabolic and hyperbolic conformal type suggests there might be a rich theory of properly embedded minimal annuli in $\nil$.

These ``horizontal catenoids" are obtained using the Weierstrass-type representation for minimal surfaces in $\nil$ \cite{dan2}. We start with a suitable harmonic map into the hyperbolic disk; this harmonic map is expressed explicitely in terms of a solution of an ODE and it will be the Gauss map of the surface. We integrate the equations and then solve a period problem. Consequently we obtain an explicit expression for these horizontal catenoids in terms of a solution of an ODE (see proposition \ref{integration}).

The second main point in this paper is to prove a half-space theorem. Discovered by Hoffman and Meeks \cite{HM}, the half-space theorem for minimal surfaces of $\R^3$ is used to understand the global geometry of proper
examples. In $\nil$, we prove that our family of ``horizontal catenoids" converges to a punctured
vertical plane and then we obtain a ``vertical half-space theorem" (vertical planes are defined in section \ref{preliminaries}).

\begin{thmintro}[theorem \ref{vhsnil}] Let $\Sigma$ be a properly immersed minimal surface in $\nil$. 
Assume that $\Sigma$ is contained on one side of a vertical plane $P$. 
Then $\Sigma$ is a vertical plane parallel to $P$.
\end{thmintro}

We next deal with complete graphs. There is a natural notion of graph in $\nil$. Indeed, $\nil$ admits a Riemannian fibration $\pi:\nil\to\R^2$ 
over the Euclidean plane. We will denote by $\xi$ a unit vector field tangent to the fibers and we will call it a vertical vector field;
it is a Killing field. Then a surface $\Sigma$ in $\nil$ is a multigraph if it is transverse to $\xi$, it is a graph if it is transverse to
  $\xi$ and $\pi_{|\Sigma}:\Sigma\to\R^2$ is injective, and it is an entire graph if it is transverse to $\xi$ and
  $\pi_{|\Sigma}:\Sigma\to\R^2$ is bijective.

A natural problem is to determine if a complete multigraph is necessarily entire. We solve this problem using our half-space theorem and applying the arguments of \cite{cr} and \cite{hrs2}. This is the following theorem.

\begin{thmintro}[theorem \ref{entire}]
Let $\Sigma$ be a complete minimal immersed surface in $\nil$. If $\Sigma$
is transverse to the vertical Killing field $\xi$, then
$\Sigma$ is an entire graph.
\end{thmintro}

Also, recently, Fernandez and Mira solved the Bernstein problem in $\nil$. We denote by $\C$ the complex plane and by $\D$ the unit disk $\{z\in\C;|z|<1\}$.

\begin{thmintro}[\cite{fermira2}]
Let $Q$ be a holomorphic quadratic differential on $\D$ or a non-identically zero holomorphic quadratic differential on $\C$. Then there exists a $2$-parameter family of generically non-congruent 
entire minimal graphs in $\nil$ whose Abresch-Rosenberg differential is $Q$. 

Conversely, all entire minimal graphs belong to these families.
\end{thmintro}

Together with our theorem \ref{entire}, this gives a classification of all complete minimal graphs in $\nil$. 

In this paper we will also deal with CMC $\frac12$ surfaces in $\h \times \R$. For these surfaces, Hauswirth, Rosenberg and Spruck \cite{hrs2} proved a half-space type theorem and used it to show that complete
multigraphs are entire, i.e., graphs over the whole hyperbolic plane $\h$. 

Our proof of theorem \ref{vhsnil} is different from the proof of the half-space type theorem in \cite{hrs2}: the main point in their proof is the construction of a continuous family of compact annuli with boundaries, contained between two horocylinders of $\h \times \R$, and converging to one of the horocylinders; they use Schauder's fixed point theorem in a quasi-linear equation; they have to control the mean curvature vector in the maximum principle. Our proof uses the family of complete annuli and the classical geometrical argument of Hoffman and Meeks \cite{HM}.

Also, by our theorem \ref{entire} and \cite{fermira2} we show that entire minimal graphs in $\nil$ correspond exactly to entire CMC $\frac12$ graphs in $\h\times\R$ by the sister surface correspondence (corollary \ref{sistergraphs}). Hence we obtain a classification of all entire CMC $\frac12$ graphs in $\h\times\R$; this solves the Bernstein problem for CMC $\frac12$ graphs in $\h\times\R$.
This is the following theorem.

\begin{thmintro}
Let $Q$ be a holomorphic quadratic differential on $\D$ or a non-identically zero holomorphic quadratic differential on $\C$. Then there exists a $2$-parameter family of generically non-congruent 
entire CMC $\frac12$ graphs in $\h \times \R$ whose Abresch-Rosenberg differential is $Q$. 

Conversely, all entire CMC $\frac12$ graphs belong to these families.
\end{thmintro}

Observe that this theorem could not be obtained using the method of \cite{fermira2}; indeed their solution of the Bernstein problem for minimal surfaces in $\nil$ is based on the relations between minimal immersions in $\nil$ and spacelike CMC immersions in Minkowski space $\LL^3$, and their arguments do not apply in our case.

We also construct a one-parameter family of properly embedded CMC $\frac12$ annuli in $\h\times\R$ which are  analogous to our minimal horizontal catenoids in $\nil$.

\begin{thmintro}[see section \ref{cmc}]
There exists a one-parameter family $(\cC_{\alpha})_{\alpha>0}$ of properly embedded CMC $\frac12$ annuli in $\h\times\R$, called ``horizontal catenoids", having the following properties:
\begin{itemize}
\item the annulus $\cC_{\alpha}$ is \emph{not} invariant by a one-parameter group of isometries,
\item the annulus $\cC_{\alpha}$ is invariant by the reflections with respect to a horizontal plane and two orthogonal vertical planes, and it is a bigraph over some domain in a horizontal plane,
\item the annulus $\cC_{\alpha}$ is conformally equivalent to $\C\setminus\{0\}$.
\end{itemize}
\end{thmintro}

The curve of intersection of $\cC_{\alpha}$ with its horizontal symmetry plane is similar to the profile curve of a rotational CMC $1$ catenoid cousin in hyperbolic space $\mathbb{H}^3$ (see \cite{bryant}, \cite{uy}). Moreover, this family converges to two punctured horocylinders tangent to each other. Hence it can give an alternative proof of the half-space type theorem of \cite{hrs2}. 

These annuli are the sister surfaces of helicoidal type minimal surfaces in $\nil$ (see section \ref{helicoid}). They are obtained in a way similar to that of the minimal horizontal catenoids in $\nil$: we start from a suitable harmonic map into the hyperbolic disk and integrate the equations of \cite{fermira2}; the period problem is solved automatically using the symmetries. Hence we obtain an explicit expression in terms of a solution of an ODE.

The paper is organized as follows. In section \ref{preliminaries} we introduce material about harmonic maps, 
minimal surfaces in $\nil$ and CMC $\frac12$ surfaces in $\h\times\R$. In section \ref{completegraphs} we give the proof and 
consequences of theorem \ref{entire}
assuming the vertical half-space theorem. In section \ref{harmonic} 
we present the family of harmonic maps that will be used in the sequel. Section \ref{horcatnil} is devoted 
to the construction of properly embedded minimal annuli in $\nil$. In section \ref{halfspace} we prove our half-space 
theorem. In section \ref{helicoid}, we construct periodic helicoidal surfaces with horizontal ``axis". In section \ref{cmc} 
we construct 
properly embedded CMC $\frac12$ annuli in $\h \times \R$.
Finally, in section \ref{appendix}
 we give the proofs of technical lemmas.

\section{Preliminaries} \label{preliminaries}

\subsection{Harmonic maps and holomorphic quadratic differentials}

In the following, we will use the unit disk
model for $\h$. We will note $\h=(\D,\si^2 (u)|\rmd u|^2)$ the disk with the hyperbolic metric 
$\si ^2 (u)|\rmd u|^2=\frac{4}{(1-|u|^2)^2}|\rmd u|^2$. The harmonic map  equation is
\begin{equation} \label{eqg}
g_{z \bar{z}}+ \frac{2\bar g}{(1-|g|^2)} g_z g_{\bar{z}}=0.
\end{equation}
\noindent
In the theory of harmonic maps there is a global
object to consider: the holomorphic quadratic
Hopf differential associated to $g$,
\begin{equation}
Q(g)= \phi (z)\rmd z ^2=(\si \circ g)^2 g_z {\bar g}_z \rmd z^{2}.
\label{eq:hopf}
\end{equation}
The function $\phi$ depends on the choice of the complex coordinate $z$, whereas $Q(g)$ does not. If $Q(g)$
is holomorphic then $g$ is harmonic. We define the function 
$\omega=\frac{1}{2} {\rm log} \frac{|g_{z}|}{|g_{\bar z}|}$.

For a given holomorphic quadratic differential $Q=\phi (z) \rmd z^2$, Wan \cite{wan} on $\D$, Wan and Au 
\cite{wan-au} on $\C$, constructed a unique (up to isometries) harmonic 
map $g :\Sigma \to \h$ with non negative Jacobian and such that the metric 
$$\tau |\rmd z|^2=4(\sigma \circ g)^2|g_{z}|^2 |\rmd z|^2=4e^{2\om} |\phi | |\rmd z|^2$$ 
is complete. To do that, they construct a spacelike CMC $\frac12$ in Minkowski space $\LL ^3$ 
with Gauss map $g$ and metric $\tau |\rmd z|^2$. First they solve the Gauss equation for the
local theory of these surfaces:
\begin{equation}
\label{eq:gordon}
\Delta_{0} \om=2 \sinh (2\om) |\phi|
\end{equation}
where $\Delta _{0}\om=4\om_{z \bar z}$. The Codazzi equation is a consequence of the fact that $\phi$ is holomorphic. 
Then a maximun principle of Cheng and Yau \cite{chengyau2} implies that there is a unique solution of (\ref{eq:gordon})
with complete
metric $\tau |\rmd z|^2$.  Then by integration of the Gauss and Codazzi equations there is a unique (up to isometries)
spacelike CMC $\frac12$ immersion 
$\tilde X=(\tilde F, \tilde h)$ in the Minkowski space $\LL ^3$. The Gauss map of $\tilde X$
is the map $g=\psi \circ \tilde N: \Sigma \to \D$,
where $\psi $ is the stereographic projection with respect to the southern pole of the quadric
$\{ |v|^2=-1\}$. The data $(Q,\tau)$ determine $g$ uniquely (up to isometries). When $\tau |\rmd z|^2$ is complete
we say that $g$ is $\tau$-complete.

In section \ref{harmonic}, we will construct a family of harmonic maps  with $Q=c\rmd z^2$ ($c \in \C$) and not
necessarily $\tau $-complete. We will use these examples to construct our horizontal
catenoids.

We describe a notion of conjugate harmonic map. It is known that a harmonic map
$g$ with $Q$ having even zeroes induces a minimal surface in $\h \times \R$. The immersion is given
$X=(g , {\rm Re}\int -2i\sqrt{Q})$ and the induced metric is 
$\rmd s^2=4\cosh^2 \omega |Q|$ (\cite{rosenhau}). Conversely, if $X=(g, t)$ is a conformal minimal immersion
then $g$ is harmonic and $Q(t)=-(t_{z})^2\rmd z^2$ is a holomorphic quadratic differential with
$Q(t)=Q(g)$. 

\begin{defin}
\label{conjugate} Two conformal minimal immersion $X,X^*:\Sigma \to \h \times \R$
are conjugate if they induce the same metric onf $\Sigma$ and if we have $Q (g^*)=-Q (g)$.
\end{defin}

In \cite{touhau} and \cite{daniel}, it is proven that the conjugate
immersion exists. If $X^*=(g^*,h^*)$, then we say that $g^*$ is the conjugate harmonic map of $g$.
In particular we will use $Q(g^*)=-Q(g)$ and $\cosh \omega ^*= \cosh \omega$ (and $\tau =\tau ^*$).

\subsection{Minimal surfaces in the Heisenberg group}
In the sequel, we  use the exponential coordinates to identify the Heisenberg group $\nil$
with $(\R^3, \rmd \sigma ^2)$, where $\rmd\sigma ^2$ given by
$$\rmd\sigma^2=\rmd x_{1}^2+\rmd x_{2}^2+\left( \rmd x_{3}  + \frac12( x_{2}\rmd x_{1}-x_{1}\rmd x_{2})\right)^2.$$

The projection $\pi:\nil\to\R^2, (x_{1},x_{2},x_{3}) \mapsto (x_{1},x_{2})$ is a Riemannian
fibration. We consider the left-invariant orthonormal frame $(E_{1},E_{2},E_{3})$ defined by
$$E_{1}=\frac{\partial}{\partial x_{1}} -\frac{x_{2}}{2}\frac{\partial}{\partial x_{3}} ,
\quad E_{2}=\frac{\partial}{\partial x_{2}}+\frac{x_{1}}{2}\frac{\partial}{\partial x_{3}} ,\quad
E_{3}= \frac{\partial}{\partial x_{3}}=\xi.$$

A vector is said to be vertical if it is proportional to $\xi$ and horizontal if it is orthogonal to $\xi$.
A surface is a multigraph if $\xi$ is nowhere tangent to it, i.e., if the restriction of
$\pi$ to the surface is a local diffeomorphism. The isometry group of $\nil$ is
$4$-dimensional and has two connected components: isometries preserving
the orientation of the fibers and the base of the fibration, and those reversing both of them.
Vertical translations are isometries. The Heisenberg group $\nil$ is a homogeneous manifold.

\begin{lemma}
\label{euclideannormal}
Let $X:\Sigma \to \nil$ be an immersion. Let $N$ be the unit normal vector to $X$
and let $\tilde N$ be the Euclidean unit normal vector to $X $ considered as an
immersion into $\R^{3}$. Then $N$ points up if and only if $\tilde N$ points up. 
\end{lemma}
\begin{proof}
We consider a conformal coordinate $z=u+iv$. In the frame $(E_{1},E_{2},E_{3})$ we have
$$X_{u}=\left[\begin{array}{c}
x_{1u} \\
x_{2u} \\
x_{3u}+\frac12 (x_{2}x_{1u}-x_{1}x_{2u})
\end{array}\right], \quad
X_{v}=\left[\begin{array}{c}
x_{1v} \\
x_{2v} \\
x_{3v}+\frac12 (x_{2}x_{1v}-x_{1}x_{2v})
\end{array}\right].
$$
Thus the third coordinate  of $X_{u} \times X_{v}$ is $x_{1u}x_{2v}-x_{1v}x_{2u}$ which is also the third
coordinate in the frame $\left(\frac{\partial}{\partial x_{1}},\frac{\partial}{\partial x_{2}},
\frac{\partial}{\partial x_{3}}\right)$ of $X_{u} \wedge X_{v}$, where $\wedge $ is the Euclidean vector product.
\end{proof}

We will call vertical planes surfaces of equation $a_1x_1+a_2x_2=b$ 
for some constants $a_1$, $a_2$ and $b$ with $(a_1,a_2)\neq(0,0)$. 
Such surfaces are minimal and flat, but not totally geodesic. Two vertical planes will be said to be parallel if their images by the projection $\pi$ are two parallel straight lines in $\R^2$.

A graph $\{x_3=f(x_1,x_2)\}$ is minimal if $f$ satisfies the 
quasi-linear equation 
$$ (1+q^2)r -2pqs+(1+p^2)t=0$$
with 
$$p=f_{x_{1}}+\frac{x_{2}}{2}, \quad q=f_{x_{2}}-\frac{x_{1}}{2},$$
$$ r=f_{x_{1}x_{1}},\quad s=f_{x_{1}x_{2}},\quad t=f_{x_{2}x_{2}}.$$

The Bernstein problem deals with the existence and the unicity of entire solutions
of this quasi-linear equation. We use conformal parametrization of  surfaces.
Let $X:\Sigma\to\nil$ be a conformal immersion. We denote by
$F=\pi\circ X$ the horizontal projection of $X$ and $h:\Sigma\to\R$
the third coordinate of $X$. We regard $F$ as a complex-valued function, identifying
$\C$ and $\R^2$. We denote the metric by $\rmd s^2=\lambda |\rmd z|^2$ and by 
$N:\Sigma\to\s^2 $ the unit normal vector
to $X$, where $\s^2$ is the unit sphere in the Lie algebra of $\nil$. 

The Gauss map of $X$ is the map $g=\psi \circ N: \Sigma \longrightarrow \bar\C=\C \cup \{\infty\}$,
where $\psi $ is the stereographic projection with respect to the southern pole, i.e.,
$g$ is defined by
$$N=\frac{1}{1+|g|^2}\left[\begin{array}{c}
2\re g \\
2\im g \\
1-|g|^2
\end{array}\right]$$
in $(E_{1},E_{2},E_{3})$. The first author proved in \cite{dan2} that the Gauss map $g$ satisfies
\begin{equation} \label{eq:harm}
(1-|g|^2)g_{z\bar z}+2\bar gg_zg_{\bar z}=0.
\end{equation}

It is important to keep in mind that $|g|=1$ exactly at points where the surface is not transverse to $\xi$.

If $\Sigma$ is a multigraph, then, up to a change of orientations, $g$ takes values in the unit disc $\D$.
When $\D$ is endowed with the hyperbolic metric $\frac{4}{(1-|z|^2)^2}|\rmd z|^2$, $g$
is a harmonic map from $\Sigma$ to $\h$. Conversely, we can recover a minimal immersion from a given harmonic
map using the following theorem.


\begin{thm}[\cite{dan2}] \label{weierstrass}
Let $\Sigma$ be a simply-connected Riemann surface.
Let $g:\Sigma\to\h$ be a harmonic map that is nowhere antiholomorphic. Let $z_0\in\Sigma$, $F_0\in\C$ and $h_0\in\R$.

Then there exists a unique conformal minimal immersion $X:\Sigma\to\nil$ such that $g$ is the Gauss map of $X$ and $X(z_0)=(F_0,h_0)$.

Moreover the immersion $X=(F,h)$ satisfies
$$F_z=-4i\frac{g_z}{(1-|g|^2)^2},\quad
F_{\bar z}=-4i\frac{g^2\bar g_{\bar z}}{(1-|g|^2)^2},$$
$$h_z=4i\frac{\bar gg_z}{(1-|g|^2)^2}-\frac i4(\bar FF_z-F\bar F_z).$$

The metric of the immersion is given by
$$\rmd s^2=16\frac{(1+|g|^2)^2}{(1-|g|^2)^4}|g_{z}|^2 | \rmd z |^2.$$
\end{thm}

The hypothesis ``nowhere antiholomorphic" forces $\lambda |\rmd z|^2$ to be a metric without
branch points. The metrics $\lambda|\rmd z|^2$ and $\tau|\rmd z|^2$ are related by 
$$\lambda=\frac{\tau}{\nu^2}$$ where $$\nu=\frac{1-|g|^2}{1+|g|^2}$$ is the third coordinate of $N$.
In the case of a multigraph we have $0<|\nu|\leqslant 1$, and so, by the above relation between $\lambda$ and $\tau$, 
it is clear that the metric $\lambda |\rmd z|^2$ is complete if $\tau |\rmd z|^2$ is complete.

In section \ref{harmonic} we will use a family of harmonic maps to construct
explicitely minimal annuli which will be the unions
of two non-complete graphs.

It is worth mentioning some recent results of Fernandez and Mira.

\begin{thm}[\cite{fermira2}] Every $\tau$-complete nowhere antiholomorphic harmonic map induces an entire minimal 
graph in $\nil$. Conversely, every entire minimal graph in $\nil$ admits a $\tau$-complete harmonic Gauss map $g$. 
\end{thm}

This theorem proves that, starting from a holomorphic quadratic differential $Q$,
there is a one-to-one canonical way to associate a two-parameter family of entire minimal graphs in $\nil$ \cite{fermira2}.

This is not enough to prove that complete multigraphs of $\nil $ are entire graphs
and then coming from a $\tau$-complete harmonic Gauss map. This fact will be the object of section \ref{completegraphs}. In other words, we will prove that, for a multigraph, if $\lambda|\rmd z|^2$ is complete, then $\tau|\rmd z|^2$ is also complete.

\subsection{Constant mean curvature $\frac12$ surfaces in $\h \times \R$}

Abresch and Rosenberg \cite{AR1} constructed a holomorphic quadratic differential $Q_{0}$ associated
to CMC $\frac12$ surfaces in $\h \times \R$; it generalizes the Hopf differential 
associated to constant mean curvature surfaces of $\R^3$. When the surface is a graph,
Fernandez and Mira \cite{fermira1} constructed a harmonic ``hyperbolic Gauss map" from the surface
to $\h$ whose associated Hopf differential is $Q=-Q_{0}$.
In addition, given a harmonic map $g$ from a surface to $\h$ plus some
additional data (described below) they construct CMC $\frac12$ graphs on $\h \times \R$
with this harmonic map as Gauss map. 

Fernandez and Mira constructed CMC $\frac12$ multigraph immersions  $X^*=(F^*,h^*):\Sigma\to
\h \times \R$ depending on the data $(Q,\tau)$.

\begin{thm}[\cite{fermira1}] \label{fermira1}
Let $\Sigma$ be a simply connected Riemann surface and $g^*:\Sigma\to\h$ be a harmonic map 
admitting data $(-Q, \tau)$. Then for any $\theta_{0}\in \C$ there exists a unique  CMC $\frac12$ immersion
$X^*=(F^*,h^*): \Sigma \to \h \times \R$ satisfying
\begin{itemize}
\item $\tau= \lambda \nu^2$, where $\lambda$ is the  conformal factor of the metric of $X^*$ and $\nu$ is the
vertical coordinate of the unit normal Gauss map,
\item $h^*_{z}(z_{0})=\theta_{0}$.
\end{itemize}
Moreover, with $G=\left( \frac{ 2 g^*}{1-|g^*|^2},\frac{1+|g^*|^2}{1-|g^*|^2}\right)$ we have
$$F^*=\frac{8{\rm Re}\left(G_{z}(4\bar Q h^*_{z} +\tau h^*_{\bar z}) \right)}{\tau ^2 -16 |Q|^2}+
 G\sqrt{\frac{\tau + 4 |h^*_{z}|^2}{\tau}}$$
and $h^*:\Sigma \to \R$ is the unique (up to an additive constant) solution to the differential
system below with $h^*_{z}(z_{0})=\theta_{0}$:
$$\left\{\begin{array}{ll}
h^*_{zz}=(\log \tau)_{z} h^*_{z} + Q\sqrt{\frac{\tau + 4 |h^*_{z}|^2}{\tau}},\\
\\
h^*_{z\bar z}=\frac{1}{4}\sqrt{\tau (\tau + 4 |h^*_{z}|^2)}.
\end{array}\right.$$
The metric can be expressed as
 $$\lambda=\frac{\tau}{\nu^2}=\tau + 4 |h^*_{z}|^2,\quad
\nu=\sqrt{\frac{\tau}{\tau + 4|h^*_{z}|^2}}.$$
\end{thm}

By the above relation between $\lambda$ and $\tau$, 
it is clear that the metric $\rmd s^2=\lambda |\rmd z|^2$ is complete if $\tau |\rmd z|^2$ is complete.
Thus, associated to a holomorphic quadratic differential $Q$,
one obtains a complete CMC $\frac12$ multigraph in $\h \times \R$. 

It is known from \cite{dan} that a CMC $\frac12$ immersion $F^*=(X^*,h^*)$
is locally isometric to a minimal immersion $X=(F,h)$ in $\nil$. These two immersions are called sister
immersions. The third coordinate $\nu$ of the unit  normal vector  of $X$ and $X^*$
remains unchanged by this correspondence. In particular the sister surface of a multigraph is a
multigraph. The harmonic Gauss maps are conjugate
($Q(g)=-Q(g^*)$ and $\tau =\tau^*$).

We mention the following result of Fernandez and Mira :

\begin{prop}[\cite{fermira2}]
If $X=(F,h)$ is a CMC  $\frac12$ minimal graph in $\h \times \R$ with a $\tau$-complete
harmonic Gauss map $g$, then $X$ is an entire graph.  
\end{prop} 

\section{Complete graphs} \label{completegraphs}
In this section we use the half-space theorem \ref{vhsnil} to obtain results on complete graphs in $\nil$ and
$\h \times \R$. 

\begin{thm} \label{entire}
Let $\Sigma$ be a complete minimal surface in $\nil$. If $\Sigma$ is transverse to the vertical Killing field $\xi$, 
then $\Sigma$ is an entire graph.
\end{thm}

\begin{cor}
Let $\Sigma$ be a complete minimal surface in $\nil$. If $\Sigma$ is transverse to the vertical Killing field $\xi$,
 then its Gauss map is $\tau$-complete.
\end{cor}

\begin{proof}
From \cite{fermira2} we know that an entire graph has a $\tau$-complete Gauss map.
\end{proof}

\begin{cor} \label{sistergraphs}
A minimal surface in $\nil$ is an entire graph if and only if
its CMC $\frac12$ sister surface in $\h \times \R$ is an entire graph.
\end{cor}

\begin{proof} By \cite{fermira2},
an entire graph of $\nil$ has a $\tau$-complete Gauss map. Then again by \cite{fermira2} the sister CMC $\frac12$
surface is entire in $\h \times \R$ (this fact comes from the completeness of $\tau | \rmd z|^2$).

Conversely, the sister of an entire CMC $\frac12$ graph in $\h \times \R$ is a complete
multigraph and then entire in $\nil$ by our theorem \ref{entire}. 
\end{proof}

\begin{cor}
Let $\Sigma$ be a complete CMC $\frac12$ surface in $\h \times \R$. If $\Sigma$ is a multigraph,
then its Gauss map is $\tau$-complete.
\end{cor}

\begin{proof}
By the theorem of Hauswirth, Rosenberg and Spruck \cite{hrs2}, $\Sigma $ is an entire graph hence
its sister also (Corollary \ref{sistergraphs}). From \cite{fermira2} we know that an entire graph is $\tau$-complete.
\end{proof}

We now prove theorem \ref{entire}. The proof is an adaptation to our case of the proof of Theorem 1.2 in \cite{hrs2}. However, we give a detailed proof for the reader's convenience and since we need to take care about the meaning of horizontal vectors in $\nil$. Lemma \ref{claim1} below is inspired by the work of Collin and Rosenberg \cite{cr}.

Let $\Sigma$ be a complete minimal surface in $\nil$ such that $\Sigma$ is transverse to the vertical Killing field $\xi$.

We assume that $\Sigma$ is not entire. We denote by $N$ the unit normal vector field to $\Sigma$. Since $\Sigma$ is a multigraph, it is orientable. The function $\nu=\langle N,\xi\rangle$ is a non-vanishing Jacobi function on $\Sigma$,
so $\Sigma$ is strongly stable and thus has bounded curvature. Hence there is $\delta >0$ such that, for each 
$p \in \Sigma$, there is a piece $G(p)$ of $\Sigma$ around $p$ that is a graph
(in exponential coordinates) over the disk $D_{2\delta}(p) \subset \rmT_{p}\Sigma$ of
radius $2\delta$ centered at the origin of $\rmT_{p}\Sigma$. This graph $G(p)$ has bounded
geometry. The $\delta$ is independant of $p$ and the bound
on the geometry of $G(p)$ is uniform as well.

We denote by $F(p)$ the image of $G(p)$ by the vertical translation mapping $p$ at height $x_3=0$, and by $F_O(p)$ the image of $G(p)$ by the translation mapping $p$ to $O=(0,0,0)\in\nil$.

In the sequel, we will call $x_3$-graphs graphs with respect to the Riemannian fibration $\pi:\nil\to\R^2$ (as explained in section \ref{preliminaries}).

We will identify vectors at different points of $\nil$ by left multiplication, and horizontal vectors with vectors in $\R^2$.

\begin{lemma} \label{claim1}
Let $(p_{n})$ be a sequence of points of $\Sigma$ such that  $N(p_n)$ has a horizontal limit $N_{\infty}$ when
$n \to +\infty$. Then there is a subsequence of $(F_O(p_{n}))_{n\in\N}$ that converges to a $\delta$-piece $P_\delta$ around $O$ of the vertical plane $P$ passing through $O$ and having $N_{\infty}$ as unit normal vector at $O$.

The convergence is in the $\cC^2$-topology. By $\delta$-piece around $O$ we mean a piece of $P$ containing all the points $p\in P$ such that $x_3(p)\in[-\delta,\delta]$ and $\pi(p)$ belongs to the closed segment of $\pi(P)\subset\R^2$ centered at $\pi(O)$ and of length $2\delta$.
\end{lemma}

\begin{proof}
Let $P$ be the vertical plane passing through $O$ and having $N_{\infty}$ as unit normal vector at $O$. We endow $P$ with the orientation induced by $N_{\infty}$.

Since the $F_O(p_n)$ have bounded geometry and are graphs over $D_{2\delta}(p_n)\subset\rmT_O(F_O(p_n))$, the $F_O(p_n)$ are bounded exponential graphs over a $\delta$-piece $P_{\delta}\subset P$ around $O$. Thus a subsequence of these graphs converges to a piece of a minimal surface $F_{\infty}$, which is tangent to $P_{\delta}$ at $O$ and which is an exponential graph over $P_{\delta}$. It suffices to show that $F_{\infty}$ is a piece of $P$.

If this is not the case, then by Theorem 5.3 in \cite{col}, in the neighbourhood
of $O$, the intersection of $F_{\infty}$ and the vertical plane consists of $m$ ($m\geqslant 2$) curves meeting at $O$. These curves separate $F_{\infty}$ into $2m$ connected components and adjacent components lie on opposite sides of the
vertical plane. Hence in a neighborhood of $O$, the Euclidean unit normal vector to $F_{\infty}$ alternates
from pointing up to pointing down as one goes from one component to the other. This is also the case for $F_O(p_{n})$, 
for $n$ large, since $F_O(p_{n})$  converges
to $F_{\infty}$ in $\cC^2$-topology. Then  by lemma \ref{euclideannormal}, 
the unit normal vector to $F_O(p_{n})$ for the metric of $\nil$ also alternates from pointing up to down. This contradicts
the fact that $F_O(p_{n})$ is transverse to $\xi$. 
\end{proof}

We now consider a piece of $\Sigma$ that is the $x_3$-graph of a function $f$ defined on the open disk $B_{R}$
of radius $R$ centered at some point $A$ of $\R^2$. Since $\Sigma$ is not an entire graph, we
choose the largest $R$ such that $f$ exists.

In the sequel, for any point $q\in B_R$ we will write $F(q)$, $N(q)$, etc. instead of $F(q,f(q))$, $N(q,f(q))$, etc.

Let $q \in \partial B_{R}$ be such that $f$ does not extend to any neighbourhood
of $q$ (to a function satisfying the minimal graph equation). 

\begin{lemma} \label{claim2}
There exists a unit horizontal vector $N_{\infty}(q)$ such that, for any sequence $q_{n} \in B_{R}$ converging to $q$, $N(q_n)\to N_{\infty}(q)$ when $n\to+\infty$. Moreover, $N_{\infty}(q)$ is normal to $\partial B_{R}$ at $q$. 
\end{lemma}

\begin{proof}
We first observe that $\nu(q_n)\to 0$ when $n\to+\infty$ (i.e. tangent planes become vertical); otherwise, the exponential graph of bounded geometry $G(q_n)$ would extend to an $x_3$-graph beyond $q$ for $q_n$ close enough to $q$, and thus the map $f$ would extend, which is a contradiction.

Let $N_{\infty}(q)$ be the horizontal unit vector at $q$ normal to $\partial B_R$ and pointing inside $B_R$. We will now prove that $N(q_n)\to N_{\infty}(q)$.

Assume that there exists a subsequence such that $N(q_n)$ converges to a horizontal vector $v\neq N_{\infty}(q)$. By lemma \ref{claim1} there exists a subsequence such that the pieces $F(q_n)$ converge to a $\delta$-piece of a vertical plane $Q$ having $v$ as unit normal vector at $q$. Since $N_{\infty}(q)$ is normal to $\partial B_{R}$ at $q$, there are points of $\pi(Q_{\delta})$ in $B_R$. Consequently there is a point $\hat q\in\pi(Q_{\delta})\cap B_R$ and a sequence $(\hat q_n)$ of points of $B_R$ converging to $\hat q$ such that $(\hat q_n,f(\hat q_n))\in G(q_n)$. Since the $F(q_n)$ converge to a $\delta$-piece of a vertical plane in the $\cC^2$-topology, $\Sigma$ has a horizontal normal at $(\hat q,f(\hat q))$, which contradicts the fact that $\Sigma$ is transverse to $\xi$.
\end{proof}

We denote by $P$ the vertical plane passing through $q$ and having $N_{\infty}(q)$ as unit normal vector at $q$.
Without loss of generality we can assume that $q=(R,0)\in\R^2$ and that $P$ is the vertical plane of equation $x_1=R$ in $\nil$. We will say that a point in $\nil$ is on the left side (respectively, on the right side) of $P$ if $x_1<R$ (respectively, $x_1>R$).

\begin{lemma} \label{claim3}
We have $f(t,0)\to\pm\infty$ when $t\to R$.
\end{lemma}

\begin{proof}
Let $\varphi(t)=f(t,0)$ and $\gamma(t)=(t,0,\varphi(t))$. 

We first claim that for $t$ close enough to $R$ we have $\varphi'(t)\neq 0$. Indeed, assume that there exists $t_0$ such that $\varphi'(t_0)=0$. We have $\gamma'(t)=E_1+\varphi'(t)\xi$ so $\gamma'(t_0)$ is horizontal. Also, $G(\gamma(t_0))$ is an exponential graph over $D_{2\delta}(\gamma(t_0))\subset\rmT_{\gamma(t_0)}\Sigma$ and $\rmT_{\gamma(t_0)}\Sigma$ contains the horizontal vector $\gamma'(t_0)$. Consequently, the projection of $G(\gamma(t_0))$ on $\R^2$ contains an open neighbourhood of $\{(t,0);t_0-\delta<t<t_0+\delta\}$. Hence, if $t_0$ is close enough to $R$, this implies that $f$ extends beyond $q=(R,0)$, which is a contradiction. This proves the claim.

Thus we can assume that $\varphi(t)$ is increasing as $t$ converges to $R$. If $\varphi(t)$ were bounded from above, then it would have a finite limit $l$ and the curve $t\mapsto(t,0,f(t,0))$ in $\Sigma$ would have finite length up till $(q,l)$. Since $\Sigma$ is complete we would have $(q,l)\in\Sigma$, but then $\Sigma$ would have a vertical tangent plane at $(q,l)$ (otherwise $f$ would extend to some neighbourhood of $q$), which gives a contradiction.
\end{proof}

From now on we assume that $f(t,0)\to+\infty$ when $t\to R$ (the case where $f(t,0)\to-\infty$ is similar).
We set $$\Gamma=\pi(P)=\{(R,s);s\in\R\}$$ and, for $\varepsilon>0$ and $s\in\R$,
$$U_{\varepsilon}=]R-\varepsilon,R[\times\R,$$
$$\gamma_{s,\varepsilon}=\{(x_1,s);R-\varepsilon<x_1<R+\varepsilon\},$$
$$\gamma^+_{s,\varepsilon}=\{(x_1,s);R-\varepsilon<x_1<R\}.$$

We fix $\varepsilon_0>0$ and consider a sequence $(t_n)$ of real numbers such that $q_n=(t_n,0)$ is in $B_R$, $q_n\to q$ when $n\to+\infty$ and
such that $$G=\bigcup_{n\in\mathbb{N}}G(q_n)$$ is connected. By lemmas \ref{claim1} and \ref{claim3} and the fact that the $G(q_n)$ are pieces of bounded geometry, $G$ is asymptotic to a part of $P$ as one goes up. Moreover we can choose the $q_n$ close enough to $R$ and to each other such that, for all $s\in[-\delta,\delta]$, the curve 
$$C_s=\pi^{-1}(\gamma_{s,\varepsilon_0})\cap G$$ is connected and has no horizontal or vertical tangents. This is possible since the $F(q_n)$ are $\cC^2$-close to $P_{\delta}$ and since $\Sigma$ is transverse to $\xi$. 



\begin{lemma} \label{claim4}
Each $G(q_n)$ is disjoint from $P$, and, for $s\in[-\delta,\delta]$, $C_s$ is an $x_3$-graph over 
$\gamma^+_{s,\varepsilon(s)}$ for some $\varepsilon(s)\in(0,\varepsilon_0]$. Moreover, $\varepsilon(s)$ can be chosed continuous.
\end{lemma}

\begin{proof}
The curve $C_s$ is an $x_3$-graph over an interval in $\gamma_{s,\varepsilon_0}$.
We show that this interval is in $\gamma^+_{s,\varepsilon_0}$. 

Suppose this is not the case for some $s_0\in[-\delta,\delta]$. Then $C_{s_0}$ has some points on the right side of $P$. But the curve $C_{0}$ stays on the left side of $P$ (otherwise $f$ would extend beyond $q$). So, for some $s_1\in]0,s_0]$, $C_{s_1}$ has points on both sides of $P$.

But $G$ is asymptotic to a part of $P$ as one goes up, so the curve $C_{s_1}$ is asymptotic to $\pi^{-1}(R,s_1)$ as the height goes to $+\infty$. This obliges $C_{s_1}$ to have a vertical tangent on the right side of $P$, which is a contradiction since $\Sigma$ is transverse to $\xi$.
\end{proof}

Consequently $\cup_{s\in[-\delta,\delta]}C_s$ is the $x_3$-graph of a function $g$ on $\cup_{s\in[-\delta,\delta]}\gamma^+_{s,\varepsilon(s)}$. The functions $f$ and $g$ coincide on the intersection of their domains of definition. The graph of $g$ on each $\gamma^+_{s,\varepsilon(s)}$ is the curve $C_s$ and the graph of $g$ is asymptotic to $P$ as the height goes to $+\infty$.


We can apply this process again replacing $C_0$ by the curve $C_{\delta}$, then $C_{-\delta}$, and so on. Analytic continuation yields an extension $h$ of $g$ to a domain $\Omega$ contained on the left side of $\Gamma$. The domain $\Omega$ is an open neighbourhood of $\Gamma$ in its left side. We have $h\to+\infty$ as one approaches $\Gamma$ in $\Omega$; the graph of $h$ is asymptotic to $P$ as the height goes to $+\infty$.

\begin{lemma} \label{claim5}
There exists $\varepsilon>0$ such that $\Omega$ contains $U_{\varepsilon}$.
\end{lemma}

\begin{proof}
The surface $\Sigma$ contains a graph over $\Omega$ composed of curves $C_s$ such that each $C_s$ is a graph over $\gamma^+_{s,\varepsilon(s)}$ for some $\varepsilon(s)>0$. Also, for each $s\in\R$, $h(t,s)$ is strictly increasing in $t$ when $t\to R$ and $h(t,s)\to+\infty$ when $t\to R$.


Let $\varepsilon_1\in(0,\delta)$ such that $\varepsilon_1\leqslant \varepsilon(s)$ for all $s\in[-\delta,\delta]$.

Suppose that for some $s_0\in\R$ we have $\varepsilon(s_0)<\varepsilon_1$. We set $\varphi(t)=h(t,s_0)$ and consider the curve $$c:(R-\varepsilon(s_0),R)\ni t\mapsto(t,s_0,\varphi(t))\in\Sigma.$$

We first claim that this curve has no horizontal tangent; indeed, if there is some $t_0$ such that $c'(t_0)$ is horizontal, then since $G(c(t_0))$ is an exponential graph over $D_{2\delta}(c(t_0))\subset\rmT_{c(t_0)}\Sigma$, and since $t_0+\delta>R$, then the curve $c$ would go on the right side of $P$, which is a contradiction.

We have $c'(t)=E_1+\left(\varphi'(t)+\frac{s_0}2\right)\xi$, so the fact that $c$ has no horizontal tangent implies that $\varphi'(t)>-\frac{s_0}2$. This implies that $\varphi(t)$ cannot tend to $+\infty$ when $t\to R-\varepsilon(s_0)$. So $\varphi(t)\to-\infty$ when $t\to R-\varepsilon(s_0)$ (otherwise the curve $c$ would extend or have a vertical tangent at $R-\varepsilon(s_0)$).

The previous discussion where we showed that the graph over $\Omega$ exists and is asymptotic to $P$ now applies to show that there is a vertical plane $\tilde P$ passing through $\tilde p=(R-\varepsilon(s_0),s_0,0)\in\nil$ such that a $\delta$-neighbourhood of $c$ in $\Sigma$ is asymptotic to a $\delta$-vertical strip in $\tilde P$ as one goes down to $-\infty$. We know this $\delta$-neighbourhood of $c$ in $\Sigma$ is asymtotic to a $\delta$-vertical strip in $P$ as one goes up to $+\infty$.

For each $s\in[s_0-\delta,s_0+\delta]$, the curve $c_s:t\mapsto(t,s,h(t,s))$ is asymptotic to some vertical line in $\tilde P$ as one goes down to $-\infty$. By analytic continuation of the $\delta$-neighbourhoods, one continues this process along $\Gamma$. 

If $P$ and $\tilde P$ are parallel, then the process continues along all of $\Gamma$ and $\Omega$ is the region bounded by $\pi(P)$ and $\pi(\tilde P)$. Then all the $\varepsilon(s)$ are equal, and this concludes the proof.

So we can assume that $P$ and $\tilde P$ intersect along some vertical line $\pi^{-1}(\hat p)$. Let us write $\hat p=(s_1,R)$. Consider the curves $c_s$ as $s$ goes from $s_0$ to $s_1$; they are graphs that become vertical both when the height goes to $+\infty$ and $-\infty$. Let $p(s)$ be the point of $C_s$ at height $0$; then, when $s\to s_1$, the path $p(s)$ has finite length (since the geometry of $\Sigma$ is bounded), so, since $\Sigma$ is complete, $p(s)$ converges to a point of $\Sigma$, and the tangent plane at this point is vertical. This contradicts the fact that $\Sigma$ is transverse to $\xi$.
\end{proof}

We can now complete the proof of theorem \ref{entire}.

\begin{proof}[Proof of theorem \ref{entire}]
We showed that $\Sigma$ contains a graph $G$ over some $U_{\varepsilon}$ which is asymptotic to $P$ as one approaches $\Gamma$ in $U_{\varepsilon}$. We apply the proof of the vertical half-space theorem \ref{vhsnil} (this theorem is stated for complete surfaces without boundary, but the proof still works in our case since $G$ is proper in some tubular neighbourhood of $P$ despite it has a non-compact boundary). This shows that such a graph $G$ cannot exist.

Consequently, $\Sigma$ is entire, and so it is an entire graph.
\end{proof}

%

\section{The family of harmonic maps}\label{harmonic}

In this section we construct a family of harmonic maps that we will use to construct annuli.
This family is derived from the two-parameter family of minimal surfaces of $\h \times \R$
constructed in \cite{hau}.

For $\alpha>0$ and $\theta\in\R$ we define $g:\C\to\bar\C$ by
$$g(u+iv)=\frac{\sin\varphi(u)+i\sinh(\alpha v+\beta(u))}{\cos\varphi(u)+\cosh(\alpha v+\beta(u))}=
\frac{\cosh (\alpha v +\beta (u)) -\cos\vp(u)}{\sin \vp(u)-i \sinh (\alpha v+\beta(u)) }$$
where $\varphi$ satisfies the following ODE:
\begin{equation} \label{eqphi}
\varphi'^2=\alpha^2+\cos(2\theta)\cos^2\varphi-\frac{\sin^2(2\theta)}{4\alpha^2}\cos^4\varphi,
\end{equation}
\noindent
and where $\beta$ is defined by
$$\beta'=\frac{\sin(2\theta)}{2\alpha}\cos^2\varphi,\quad \beta(0)=0.$$
The function $\vp$ is defined on the whole $\R$. We will study this function $\varphi$ in lemma \ref{lemmaphi}.
We also set 
$$A=\alpha v+\beta(u),\quad D=\cos\varphi+\cosh A.$$ 
We notice that
\begin{equation} \label{modg}
1-|g|^2=\frac{2\cos\varphi}D.
\end{equation}

\begin{prop} \label{hopf}
The function $g$ satisfies 
$$(1-|g|^2)g_{z \bar z} + 2 \bar g g_{z} g_{\bar z}=0$$
and its Hopf differential is
$$Q=\frac14e^{-2i\theta}\rmd z^2.$$
\end{prop}

\begin{proof}
To see that $g$ satisfies the equation, it suffices to see that 
$$Q=\frac4{(1-|g|^2)^2}g_{z}g_{\bar z}\rmd z^2$$
is holomorphic.

We compute 
$$g_u=\frac{\varphi'+i\beta'}{D^2}(1+\cos\varphi\cosh A+i\sin\varphi\sinh A),$$
$$g_v=\frac{i\alpha}{D^2}(1+\cos\varphi\cosh A+i\sin\varphi\sinh A).$$
From this and \eqref{modg} we get 
$$Q=\frac{{\varphi'}^2-(\alpha+i\beta')^2}{4\cos^2\varphi}\rmd z^2.$$
Using \eqref{eqphi} and the definition of $\beta$ we get $Q=\frac14e^{-2i\theta}\rmd z^2$.
\end{proof}

For $\alpha>0$ and $\theta\in\R$, we set
$$C=C_{\alpha,\theta}=\frac{\sin(2\theta)}{2\alpha},\quad
P_{\alpha,\theta}(x)=\alpha^2+\cos(2\theta)x^2-C_{\alpha,\theta}^2x^4,$$
so that \eqref{eqphi} is equivalent to $${\varphi'}^2=P_{\alpha,\theta}(\cos\varphi).$$

We set $\theta^+_{\alpha}=\frac\pi2$ if $\alpha>1$ and
$\theta^+_{\alpha}=\frac12\arccos(1-2\alpha^2)\in(0,\frac\pi2]$ if $\alpha\leqslant 1$. Let $\Omega=\{(\alpha,\theta)\in\R^2;\alpha>0,\theta\in(-\theta^+_{\alpha},\theta^+_{\alpha})\}$.

If $2\theta\notin\pi\Z$, we have
$$P_{\alpha,\theta}(x)=C_{\alpha,\theta}^2(\rho^-_{\alpha,\theta}-x^2)(\rho^+_{\alpha,\theta}+x^2)$$
with $$\rho^-_{\alpha,\theta}=\frac{2\alpha^2}{1-\cos(2\theta)},\quad
\rho^+_{\alpha,\theta}=\frac{2\alpha^2}{1+\cos(2\theta)}.$$ Thus, if $2\theta\notin\pi\Z$ and $(\alpha,\theta)\in\Omega$, then $\rho^-_{\alpha,\theta}>1$. Also, we have $P_{\alpha,0}(x)=\alpha^2+x^2$. From this we deduce that 
$$\forall(\alpha,\theta)\in\Omega,\forall x\in[-1,1],P_{\alpha,\theta}(x)>0.$$
Thus, if $(\alpha,\theta)\in\Omega$, then the right term in \eqref{eqphi} does not vanish.

\begin{lemma} \label{lemmaphi}
Let $(\alpha,\theta)\in\Omega$. Let $\varphi$ be the solution of \eqref{eqphi} such that $\varphi(0)=0$ and $\varphi'(0)\leqslant 0$. Then
\begin{itemize}
\item[1.] $\forall u,\varphi'(u)<0$,
\item[2.] the function $\varphi$ is a decreasing bijection from $\R$ onto $\R$,
\item[3.] there exists a real number $U>0$ such that
$$\forall u\in\R,\varphi(u+U)=\varphi(u)-\pi,$$
\item[4.] the function $\varphi$ is odd.
\end{itemize}
\end{lemma}

\begin{proof}
\begin{itemize}
\item[1.] Since the right term in \eqref{eqphi} does not vanish, $\varphi'$ does not vanish.
\item[2.] The right term in \eqref{eqphi} is bounded by two positive constants $c_1$ and $c_2$; hence $-\sqrt{c_2}\leqslant\varphi'\leqslant-\sqrt{c_1}$, which proves that $\varphi$ is defined on the entire $\R$ and that $\varphi(u)\to-\infty$ when $u\to+\infty$ and $\varphi(u)\to+\infty$ when $u\to-\infty$.
\item[3.] There exists $U>0$ such that $\varphi(U)=-\pi$. Then the function $\tilde\varphi:u\mapsto\varphi(u+U)+\pi$ satisfies \eqref{eqphi} with $\tilde\varphi(0)=0$ and $\tilde\varphi'<0$; hence $\tilde\varphi=\varphi$.
\item[4.] The function $\hat\varphi:u\mapsto-\varphi(-u)$ satisfies \eqref{eqphi} with $\hat\varphi(0)=0$ and $\hat\varphi'<0$; hence $\hat\varphi=\varphi$.
\end{itemize}
\end{proof}

In the sequel we will use the function $G:\R\to\R$ defined by
$$G'=\frac{C^2\cos^2\varphi-\cos(2\theta)}{\alpha-\varphi'},\quad G(0)=0.$$
(We recall that $\alpha-\varphi'>0$.) The functions $\beta$ and $G$ are odd and satisfy
\begin{equation} \label{periodbetaG}
\forall u\in\R,\beta(u+U)=\beta(u)+\beta(U),\quad
\forall u\in\R,G(u+U)=G(u)+G(U).
\end{equation}

\begin{lemma}
We have 
\begin{equation} \label{U2}
\varphi\left(\frac U2\right)=-\frac{\pi}2,\quad
\beta\left(\frac U2\right)=\frac{\beta(U)}2,\quad
G\left(\frac U2\right)=\frac{G(U)}2.
\end{equation}
\end{lemma}

\begin{proof}
We have $\varphi\left(\frac U2\right)=\varphi\left(-\frac U2\right)-\pi=-\varphi\left(\frac U2\right)-\pi$, which gives the first formula. We prove the other formulas in the same way.
\end{proof}

\begin{lemma}
The following identities hold.
\begin{equation} \label{phiprimealpha}
\varphi'+\alpha=G'\cos^2\varphi,
\end{equation}
\begin{equation} \label{phisecond}
\varphi''=-(\cos(2\theta)-2C^2\cos^2\varphi)\sin\varphi\cos\varphi,
\end{equation}
\begin{equation} \label{Gsecondcosphi}
G''\cos\varphi=(2\varphi'G'-\cos(2\theta)+2C^2\cos^2\varphi)\sin\varphi,
\end{equation}
\begin{equation} \label{Gsecond}
G''=\frac{2C^2\alpha-\cos(2\theta)G'}{\alpha-\varphi'}\sin\varphi\cos\varphi,
\end{equation}
\begin{equation} \label{Gsecond2}
G''=(C^2+{G'}^2)\sin\varphi\cos\varphi.
\end{equation}
\end{lemma}

\begin{proof}
Formulas \eqref{phiprimealpha}, \eqref{phisecond}, \eqref{Gsecondcosphi} and \eqref{Gsecond} are straightforward. 
Using \eqref{Gsecond} and the definition of $G$ we get
$$(\varphi'-\alpha)^2G''=(\cos^2(2\theta)-C^2\cos(2\theta)\cos^2\varphi-2\alpha C^2(\varphi'-\alpha))
\sin\varphi\cos\varphi.$$
On the other hand we have
$$(\varphi'-\alpha)^2(C^2+{G'}^2)=\cos^2(2\theta)-C^2\cos(2\theta)\cos^2\varphi-2\alpha C^2(\varphi'-\alpha).$$
This proves \eqref{Gsecond2}.
\end{proof}

For $(\alpha,\theta)\in\Omega$ we set
$$L(\alpha,\theta)=\int_{-1}^1
\frac{2\alpha C_{\alpha,\theta}^2x^2-\alpha\cos(2\theta)+C_{\alpha,\theta}^2x^2\sqrt{P_{\alpha,\theta}(x)}}
{\sqrt{(1-x^2)P_{\alpha,\theta}(x)}(\alpha+\sqrt{P_{\alpha,\theta}(x)})}\rmd x.$$
We will prove in section \ref{appendix} the following technical lemmas.

\begin{lemma} \label{lemmaL}
Let $\alpha>0$. Then there exists a unique $\tilde\theta_{\alpha}\in(0,\theta^+_{\alpha})\cap(0,\frac\pi4)$ such that 
$$L(\alpha,\tilde\theta_{\alpha})=0.$$
\end{lemma}

\begin{lemma} \label{limitthetainfinity}
We have $$\lim_{\alpha\to+\infty}\tilde\theta_{\alpha}=\frac{\pi}4.$$
\end{lemma}

%
%
%
%

\section{Horizontal catenoids in $\nil$} \label{horcatnil}

In this section we construct a one-parameter  family of properly embedded
minimal annuli in $\nil$. We use the notations of section \ref{harmonic}.

We will start from the map $g$ which satisfies \eqref{eqg} (by proposition \ref{hopf}) outside points where $|g|=1$ but which
 does not take values in $\D$. However, in this case we can still recover a minimal immersion (but not a multigraph) by 
 theorem \ref{weierstrass} provided the map we obtain is well-defined when $|g|=1$ and provided the metric we obtain 
 has no singularity. In fact these two kinds of problems do not appear in our case, as shown by the following proposition.

\begin{prop}
\label{integration}
The conformal minimal immersion $X=(F,h):\C\to\nil$ whose Gauss map is $g$ is given (up to a translation) by
$$F(u+iv)=\frac{G'}{\alpha}\cos\varphi\sinh A-\frac{C}{\alpha}\sin\varphi\cosh A
+i(Cv-G),$$
\begin{eqnarray*}
h(u+iv) & = & -\frac1{\alpha}\left(G'\sin\varphi+\frac{C^2}{\alpha}\sin\varphi
+\frac{(Cv-G)G'}2\cos\varphi\right)\sinh A \\
& & +\frac1{\alpha}\left(-C\cos\varphi+\frac{CG'}{\alpha}\cos\varphi
+\frac{C(Cv-G)}2\sin\varphi\right)\cosh A.
\end{eqnarray*}

The metric of $X$ is given by 
$$\rmd s^2=({G'}^2+C^2)\cosh^2A|\rmd z|^2.$$
\end{prop}

\begin{proof}
We first recover $F$ using theorem \ref{weierstrass} and the above computations. We get
$$F_z=-\frac i{2\cos^2\varphi}(\varphi'+i\beta'+\alpha)(1+\cos\varphi\cosh A+i\sin\varphi\sinh A),$$
$$F_{\bar z}=-\frac i{2\cos^2\varphi}(\varphi'-i\beta'+\alpha)(1-\cos\varphi\cosh A+i\sin\varphi\sinh A),$$
hence
\begin{eqnarray*}
F_u & = & \frac{\beta'\cos\varphi\cosh A-i(\varphi'+\alpha)(1+i\sin\varphi\sinh A)}{\cos^2\varphi} \\
& = & C\cos\varphi\cosh A-iG'(1+i\sin\varphi\sinh A),
\end{eqnarray*}
\begin{eqnarray*}
F_v & = & \frac{(\varphi'+\alpha)\cos\varphi\cosh A+i\beta'(1+i\sin\varphi\sinh A)}{\cos^2\varphi} \\
& = & G'\cos\varphi\cosh A+iC(1+i\sin\varphi\sinh A).
\end{eqnarray*}
This gives $F$.

Then we get
\begin{eqnarray*}
h_z & = & \frac{G'+iC}4(2\cos\varphi\sinh A+2i\sin\varphi\cosh A \\
& & -\frac{G'}{\alpha}\cos\varphi\sinh A+\frac C{\alpha}\sin\varphi\cosh A \\
& & +i(Cv-G)\cos\varphi\cosh A-(Cv-G)\sin\varphi\sinh A).
\end{eqnarray*}
This gives $h$.

Using \eqref{modg} and computations done in the proof of proposition \ref{hopf} we get
$$1+|g|^2=\frac{2\cosh A}D,\quad |g_z|^2=\frac{(\varphi'+\alpha)^2+{\beta'}^2}{4D^2},$$
and so by theorem \ref{weierstrass} we obtain the formula.
\end{proof}

\begin{prop} \label{period}
Let $\alpha>0$ and $\theta=\tilde\theta_{\alpha}$. Then the corresponding immersion $X$ is simply periodic, i.e., there exists $Z\in\C\setminus\{0\}$ such that
$$\forall z\in\C,X(z+Z)=X(z).$$
\end{prop}

\begin{proof}
Let $C_{\alpha}=C_{\alpha,\tilde\theta_{\alpha}}$ and $P_{\alpha}(x)=P_{\alpha,\tilde\theta_{\alpha}}(x)$. We set $$V=-\frac{\beta(U)}{\alpha}.$$ Then, by \eqref{periodbetaG}, for all $(u,v)\in\R^2$, we have $A(u+U+i(v+V))=A(u+iv)$.

We claim that, for all $(u,v)\in\R^2$, we have $\im F(u+U+i(v+V))=\im F(u+iv)$, i.e., that
\begin{equation} \label{period0}
\alpha G(U)+C\beta(U)=0.
\end{equation}
We have $$G(U)=\int_0^UG'(u)\rmd u,\quad
\beta(U)=\int_0^U\beta'(u)\rmd u.$$
We do the change of variables $x=\cos\varphi(u)$, hence $\rmd x=-\varphi'\sin\varphi\rmd u=\varphi'\sqrt{1-x^2}\rmd u$ since $\varphi\in[-\pi,0]$. We get
$$G(U)=\int_{-1}^1\frac{C_{\alpha}^2x^2-\cos(2\tilde\theta_{\alpha})}
{\sqrt{(1-x^2)P_{\alpha}(x)}(\alpha+\sqrt{P_{\alpha}(x)})}\rmd x,$$
$$\beta(U)=\int_{-1}^1\frac{C_{\alpha}x^2}
{\sqrt{(1-x^2)P_{\alpha}(x)}}\rmd x,$$
and so $\alpha G(U)+C\beta(U)=L(\alpha,\tilde\theta_{\alpha})=0$  by lemma \ref{lemmaL}. This proves the claim.

Hence $A(u+iv)$ and $\im F(u+iv)=Cv-G(u)$ are $(U+iV)$-periodic. We set $Z=2(U+iV)$ (we have $Z\neq 0$ since $U>0$). Then it 
follows from the expressions of $F$ and $h$ that they are $Z$-periodic.
\end{proof}

\begin{defn}
Let $\alpha>0$. The surface given by $X$ when $\theta=\tilde\theta_{\alpha}$ is called a horizontal catenoid of parameter $\alpha$ with respect to the $x_2$-axis. It will be denoted $\cC_{\alpha}$.
\end{defn}

The coordinates $(x_1,x_2,x_3)$ of $\cC_{\alpha}$ are
\begin{eqnarray*}
x_1 & = & \frac{G'(u)}{\alpha}\cos\varphi(u)\sinh A-\frac C{\alpha}\sin\varphi(u)\cosh A, \\
x_2 & = & \frac{C}{\alpha}A-\frac{C}{\alpha}\beta(u)-G(u), \\
x_3 & = & -\frac{x_1x_2}2+\frac C{\alpha}\left(\frac{G'(u)}{\alpha}-1\right)\cos\varphi(u)\cosh A \\
& & -\frac1{\alpha}\left(\frac{C^2}{\alpha}+G'(u)\right)\sin\varphi(u)\sinh A.
\end{eqnarray*}

We now study the geometry of $\cC_{\alpha}$.
We first notice that $$\left\{\begin{array}{lll}
x_1(u+U,v+V) & = & -x_1(u,v), \\
x_2(u+U,v+V) & = & x_2(u,v), \\
x_3(u+U,v+V) & = & -x_3(u,v),
\end{array}\right.$$ 
so $\cC_{\alpha}$ is invariant by the rotation of angle $\pi$ around the $x_2$-axis.
We also have
$$\left\{\begin{array}{lll}
x_1(-u,-v) & = & -x_1(u,v), \\
x_2(-u,-v) & = & -x_2(u,v), \\
x_3(-u,-v) & = & x_3(u,v),
\end{array}\right.$$  
so $\cC_{\alpha}$ is invariant by the rotation of angle $\pi$ around the $x_3$-axis.
Since the composition of the rotations of anngle $\pi$ around the $x_2$ and $x_3$ axes is the rotation of angle $\pi$ around the $x_1$-axis, $\cC_{\alpha}$ is also invariant by this rotation.


It will be convenient to use the following coordinates in $\nil$:
\begin{equation}
\label{eq:coordinates}
y_1=x_1,\quad y_2=x_2,\quad y_3=x_3+\frac{x_1x_2}2.
\end{equation}
In these coordinates the metric of $\nil$ is given by
$$\rmd y_1^2+\rmd y_2^2+(\rmd y_3-y_1\rmd y_2)^2.$$ 
In particular, in a vertical plane of equation $y_2=c$ ($c\in\R$), the pair $(y_1,y_3)$ is a pair of Euclidean coordinates.

We now study the intersection of $\cC_{\alpha}$ with a vertical plane of equation $y_2=c$ ($c\in\R$). On $\cC_{\alpha}$, this intersection is given by 
\begin{equation} \label{eqA}
A=\frac{\alpha}Cc+\beta(u)+\frac{\alpha}CG(u).
\end{equation}
Hence, reporting this equality in the expressions of $(x_1,x_2,x_3)$, we obtain a parametrization $u\mapsto\gamma(u)$ of this intersection.

\begin{lemma} \label{x2constant}
On a curve where $y_2$ is constant we have
$$y_1'(u)=\frac{C^2+{G'}^2}{C}\cos\varphi\cosh A,$$
$$y_3'(u)=-\frac{C^2+{G'}^2}{C}\sin\varphi\cosh A.$$
\end{lemma}

\begin{proof}
Differentiating \eqref{eqA} we obtain $A'=C\cos^2\varphi+\frac{\alpha}CG'$. Hence we get
\begin{eqnarray*}
y_1'(u) & = & \frac1{\alpha}\left(G''\cos\varphi-G'\varphi'\sin\varphi-C^2\sin\varphi\cos^2\varphi
-\alpha G'\sin\varphi\right)\sinh A \\
& & +\frac1{\alpha}\left(CG'\cos^3\varphi+\frac{\alpha}C{G'}^2\cos\varphi-C\varphi'\cos\varphi\right)
\cosh A
\end{eqnarray*}
and
\begin{eqnarray*}
y_3'(u) & = & \frac1{\alpha}\left(\frac C{\alpha}G''\cos\varphi-C\left(\frac{G'}{\alpha}-1\right)\varphi'\sin\varphi \right.\\
& &\left. -\left(\frac{C^2}{\alpha}+G'\right)\left(C\cos^2\varphi+\frac{\alpha}CG'\right)\sin\varphi\right)\cosh A \\
& & +\frac1{\alpha}\left(C\left(\frac{G'}{\alpha}-1\right)\left(C\cos^2\varphi+\frac{\alpha}CG'\right)\cos\varphi
\right.\\
& &\left. -G''\sin\varphi-\left(\frac{C^2}{\alpha}+G'\right)\varphi'\cos\varphi\right)\sinh A.
\end{eqnarray*}
We conclude using \eqref{phiprimealpha} and \eqref{Gsecond2}.
\end{proof}

\begin{prop} \label{convex}
Let $c\in\R$. The intersection of $\cC_{\alpha}$ and the vertical plane $\{y_2=c\}$ is a non-empty closed embedded convex curve.
\end{prop}

\begin{proof}
This intersection is non-empty since setting $u=0$ and $A=c$ gives $y_2=c$. Also, by lemma \ref{x2constant} we have ${y_1'}^2+{y_3'}^2>0$, so the intersection of $\cC_{\alpha}$ and the vertical plane $\{y_2=c\}$ is a smooth curve $\gamma$. Also, we have $\gamma(u+2U)=\gamma(u)$, so the curve is closed.

We now prove that $\gamma$ is embedded and convex. We consider the half of $\gamma$ corresponding to $u\in(-\frac U2,\frac U2)$. We have $\cos\varphi(u)>0$. Then, by lemma \ref{x2constant}, $u\mapsto y_1(u)$ is injective and increasing. We get
$$\frac{\rmd y_3}{\rmd y_1}=-\tan\varphi(u),$$ so $\frac{\rmd y_3}{\rmd y_1}$ is an increasing function of $u$, and also of $y_1$. Consequently, the half of $\gamma$ corresponding to $u\in(-\frac U2,\frac U2)$ is an embedded convex arc and is situated below the segment linking its endpoints.

Finally, since $\gamma(u+U)=-\gamma(u)$, the whole curve is embedded and convex.
\end{proof}



\begin{thm} \label{propertieshorcat}
The horizontal catenoid $\cC_{\alpha}$ has the following properties.
\begin{itemize}
\item[1.] The intesection of $\cC_{\alpha}$ and any vertical plane of equation $x_2=c$ ($c\in\R$) is a non-empty closed embedded convex curve.
\item[2.] The surface $\cC_{\alpha}$ is properly embedded.
\item[3.] The horizontal catenoid $\cC_{\alpha}$ is invariant by rotations of angle $\pi$ around the $x_1$, $x_2$ and $x_3$ axes. The $x_2$-axis is contained in the ``interior'' of $\cC_{\alpha}$.
\item[4.] It is conformally equivalent to $\C\setminus\{0\}$.
\end{itemize}
\end{thm}

\begin{proof}
\begin{itemize}
\item[1.] This is proposition \ref{convex}.
\item[2.] The fact that $\cC_{\alpha}$ is embedded is a consequence of proposition \ref{convex}. On a diverging path on $\cC_{\alpha}$, $A$ must be diverging and so $x_2$ is diverging. Consequently, $\cC_{\alpha}$ is proper.
\item[3.] The symmetries of $\cC_{\alpha}$ have already been proved.
The $x_2$-axis is contained in the ``interior'' of $\cC_{\alpha}$ since each curve $x_2=c$ ($c\in\R$) is convex and symmetric with respect to the $x_2$-axis.
\item[4.] The immersion $X=(F,h)$ induces a conformal bijective parametrization of $\cC_{\alpha}$ by $\C/(\Z Z)$.
\end{itemize}
\end{proof}

We now discribe a few remarkable curves on $\cC_{\alpha}$.

The curve corresponding to $u=0$ is the set of to the lowest points of the curves $y_2=c$ ($c\in\R$). This curve is given by
$$\left\{\begin{array}{lll}
y_1 & = & \frac{\alpha-\sqrt{\alpha^2+\cos(2\theta)-C^2}}{\alpha}\sinh\left(\frac{\alpha}Cy_2\right), \\
y_3 & = & -\frac{C\sqrt{\alpha^2+\cos(2\theta)-C^2}}{\alpha}\cosh\left(\frac{\alpha}Cy_2\right).
\end{array}\right.$$

The curves along which $\cC_{\alpha}$ is vertical correspond to $u=\pm\frac U2$ (because of formula \eqref{modg}). They are symmetric one to the other with respect to the $x_2$-axis. By \eqref{period0} and \eqref{U2}, the curve corresponding to $u=\frac U2$ is given by
$$\left\{\begin{array}{lll}
y_1 & = & \frac{C}{\alpha}\cosh\left(\frac{\alpha}Cy_2\right), \\
y_3 & = & \frac{2C^2-\cos(2\theta)}{2\alpha^2}\sinh\left(\frac{\alpha}Cy_2\right).
\end{array}\right.$$
Consequently, the horizontal projection of $\cC_{\alpha}$ is
$$\pi(\cC_{\alpha})
=\left\{(y_1,y_2)\in\R^2;|y_1|\leqslant\frac C{\alpha}\cosh\left(\frac{\alpha}Cy_2\right)\right\}$$ It is a remarkable fact that this projection coincides with the projection of a minimal catenoid of $\R^3$ of parameter $\frac C{\alpha}$.

The curve given by $x_2=0$ is the analog of the ``waist circle'' of minimal catenoids in $\R^3$.

\begin{prop}
On $\cC_{\alpha}$, there exists some points with negative curvature and some points with positive curvature. 
Moreover, $\cC_{\alpha}$ has infinite total absolute curvature.
\end{prop}

\begin{proof}
Setting $\lambda=({G'}^2+C^2)\cosh^2A$, the curvature of $\rmd s^2$ is given by
$$K=-\frac1{2\lambda}\Delta_0(\ln\lambda)$$ where $\Delta_0$ is the Laplacian with respect to $|\rmd z|^2$. Thus we have
\begin{eqnarray*}
K\lambda & = & -\frac{\partial}{\partial u}(\beta'\tanh A)-\frac{\partial}{\partial v}(\alpha\tanh A)
-\frac{\partial}{\partial u}\left(\frac{G'G''}{{G'}^2+C^2}\right) \\
& = & 2C\varphi'\sin\varphi\cos\varphi\tanh A-\frac{C^2\cos^4\varphi+\alpha^2}{\cosh^2A} \\
& & -(C^2+{G'}^2)\sin^2\varphi\cos^2\varphi-G'\varphi'(2\cos^2\varphi-1)
\end{eqnarray*}
by \eqref{Gsecond2}.

Hence, when $u=\pm\frac U2$ we have $K\lambda=-\frac{\alpha^2}{\cosh^2A}+\frac{\cos(2\theta)}2$, which is positive for $|A|$ large enough. On the other hand, when $u=A=0$ we get $K\lambda=-2\alpha^2-\cos(2\theta)+\alpha\sqrt{\alpha^2+\cos(2\theta)-C^2}<0$.

Finally, the total absolute curvature of $\cC_{\alpha}$ is
$$\int_{-U}^U\int_{-\infty}^{+\infty}|K|\lambda\rmd u\rmd v=+\infty$$ since, in general,
$K\lambda$ does not tend to $0$ when $v\to+\infty$ and $u$ fixed.
\end{proof}

\section{Limit of horizontal catenoids and vertical half-space theorem in $\nil$}
\label{halfspace}

%
%

%

In this section we study the limit of $\cC_{\alpha}$ when $\alpha\to+\infty$. As a corollary we obtain a vertical half-space theorem.

Since the parameter $\alpha$ will vary, the quantities and functions appearing in the construction of $\cC_{\alpha}$ will 
be denoted by $X_{\alpha}$, $\varphi_{\alpha}$, $U_{\alpha}$, $C_{\alpha}$, $\beta_{\alpha}$, etc. instead of $X$, $\varphi$, $U$, $C$,
 $\beta$, etc. They depend smoothly on $\alpha$. 

\begin{prop} \label{limitinfinity}
Let $(\hat u,\hat v)\in\mathbb{R}^2$. For $\alpha>0$, let $u_{\alpha}=\frac{\hat u}{\alpha}$ and $v_{\alpha}=\frac{4\ln\alpha+\hat v}{\alpha}$. Then, when $\alpha\to+\infty$,
$$(y_1)_{\alpha}(u_{\alpha},v_{\alpha})\to\frac{\sin\hat u}4e^{\hat v/2},$$
$$(y_2)_{\alpha}(u_{\alpha},v_{\alpha})\to0,$$
$$(y_3)_{\alpha}(u_{\alpha},v_{\alpha})\to-\frac{\cos\hat u}4e^{\hat v/2}.$$
\end{prop}

\begin{proof}
We have $\tilde\theta_{\alpha}\to\frac{\pi}4$ and so $C_{\alpha}\sim\frac1{2\alpha}$.
We have $|(\beta_{\alpha})'|\leqslant\frac{1}{2\alpha}$, thus $|\beta_{\alpha}(u)|\leqslant\frac{|u|}{2\alpha}$ and so $\beta_{\alpha}(u_{\alpha})=\rmO\left(\frac1{\alpha^2}\right)$. Also, for $\alpha\geqslant \frac12$, we have 
$$|(G_{\alpha})'|=
\left|\frac{C_{\alpha}^2\cos^2\varphi_{\alpha}-\cos(2\tilde\theta_{\alpha})}{\alpha-(\varphi_{\alpha})'}\right|
\leqslant\frac{1}{4\alpha^3}+\frac{\cos(2\tilde\theta_{\alpha})}{\alpha},$$
thus $G_{\alpha}(u_{\alpha})=\rmO\left(\frac1{\alpha^2}\right)$. From this we obtain that $$(y_2)_{\alpha}(u_{\alpha},v_{\alpha})=C_{\alpha}v_{\alpha}-G_{\alpha}(u_{\alpha})\to 0.$$

We also have $$A=4\ln\alpha+\hat v+\rmO\left(\frac1{\alpha^2}\right),\quad 
\cosh A\sim\frac12e^{\hat v/2}\alpha^2,\quad \sinh A\sim\frac12e^{\hat v/2}\alpha^2.$$
And, since, for $\alpha\geqslant 1$, $-\sqrt{\alpha^2+1}\leqslant(\varphi_{\alpha})'\leqslant-\sqrt{\alpha^2-1}$, we have $\varphi_{\alpha}(u_{\alpha})\to-\hat u$. This concludes the proof.
\end{proof}

This propostition means that, when $\alpha\to+\infty$, the half of $\cC_{\alpha}$ corresponding to $A>0$ converges to the punctured 
vertical plane $\{x_2=0\}\setminus\{(0,0,0)\}$. In the same way one can prove that the other half of $\cC_{\alpha}$ converges to this 
punctured vertical plane. 

\begin{lemma} \label{waistcircle}
The curve of equation $y_2=0$ in $\cC_{\alpha}$ converges uniformly to $0$ when $\alpha\to+\infty$.
\end{lemma}

\begin{proof}
On the curve of equation $y_2=0$ in $\cC_{\alpha}$ we have
$$|(y_1)_{\alpha}(u)|\leqslant\frac{|G_{\alpha}'(u)|}{\alpha}\sinh|A_{\alpha}(u)|
+\frac{C_{\alpha}}{\alpha}\cosh A_{\alpha}(u),$$
$$|(y_3)_{\alpha}(u)|\leqslant
\frac{C_{\alpha}}{\alpha}\left(\frac{|G_{\alpha}'(u)|}{\alpha}+1\right)\cosh A_{\alpha}(u)
+\frac1{\alpha}\left(\frac{C_{\alpha}^2}{\alpha}+|G_{\alpha}'(u)|\right)\sinh|A_{\alpha}(u)|,$$
with $$A_{\alpha}(u)=\beta_{\alpha}(u)+\frac{\alpha}{C_{\alpha}}G_{\alpha}(u).$$
By the computations done in the proof of proposition \ref{limitinfinity} we have, for $\alpha\geqslant 1$. $|G_{\alpha}'(u)|\leqslant\frac2\alpha$ and $C_{\alpha}\leqslant\frac1{2\alpha}$. Hence it suffices to prove that $A_{\alpha}$ is uniformly bounded when $\alpha\to+\infty$.

By \eqref{period0}, the function $\beta_{\alpha}+\frac{\alpha}{C_{\alpha}}G_{\alpha}$ is $2U_\alpha$-periodic with 
$$U_{\alpha}=\int_0^{U_{\alpha}}\rmd u=\int_{-1}^1\frac{\rmd x}{\sqrt{(1-x^2)P_{\alpha}(x)}}.$$
We now assume that $\alpha\leqslant 1$. For $x\in[-1,1]$ we have $P_{\alpha}(x)\geqslant\alpha^2-1$, and so $U_{\alpha}\leqslant\frac{\pi}{\sqrt{\alpha^2-1}}$. Using the bounds on $\beta_{\alpha}'$ and $G_{\alpha}'$ and the periodicity we get $$|A_{\alpha}(u)|\leqslant\frac{\pi}{2\alpha\sqrt{\alpha^2-1}}
+\frac{\alpha}{C_{\alpha}}\frac{2\pi}{\alpha\sqrt{\alpha^2-1}}.$$
Next, since $C_{\alpha}\sim\frac1{2\alpha}$, we conclude that $A_{\alpha}$ is uniformly bounded when $\alpha\to+\infty$, which ends the proof.
\end{proof}

\begin{thm}[vertical half-space theorem] \label{vhsnil}
Let $\Sigma$ be a properly immersed minimal surface in $\nil$. 
Assume that $\Sigma$ is contained on the one side of a vertical plane $P$. 
Then $\Sigma$ is a vertical plane parallel to $P$.
\end{thm}

\begin{proof}
We assume that $\Sigma$ is not a vertical plane.

We proceed as in \cite{HM}. Up to an isometry of $\nil$ we can assume that $P$ is the plane 
$\{y_2=0\}$, that $\Sigma\subset\{y_2\leqslant 0\}$ and that $\Sigma$ is not contained in any half-space 
$\{y_2\leqslant-\varepsilon\}$ for $\varepsilon>0$. By the maximum principle, we necessarily have $\Sigma \cap P=\emptyset$.

We use the coordinates $(y_{1},y_{2},y_{3})$ defined by \eqref{eq:coordinates}.
For $\varepsilon\in\R$, let $T_{\varepsilon}:(y_1,y_2,y_3)\mapsto(y_1,y_2+\varepsilon,y_3)$ 
(this is a translation in the $y_2$ direction, an isometry of $\nil$). 
Then, for $\varepsilon>0$ sufficiently small, we have $T_{\varepsilon}(\Sigma)\cap P\neq\emptyset$.

For $\alpha\geqslant 1$ we consider the half-horizontal catenoid $\cC'_{\alpha}=\cC_{\alpha}\cap\{y_2\geqslant 0\}$. 
By lemma \ref{waistcircle}, there exists a compact subset $\cD$ of $P$ containing $0$ and $\cC_{\alpha}\cap P$ for all $\alpha\geqslant 1$.

We claim that there exists $\varepsilon>0$ such that
$$T_{\varepsilon}(\Sigma)\cap P\neq\emptyset,\quad
T_{\varepsilon}(\Sigma)\cap\cC'_1=\emptyset,\quad
T_{\varepsilon}(\Sigma)\cap\cD=\emptyset.$$ 
Assume the claim is false. Since $T_{\eta}(\Sigma)\cap P\neq\emptyset$ for $\eta$ small enough, this means that 
there exists a sequence $(\varepsilon_n)$ of positive numbers converging to $0$ and a sequence $(q_n)$ of points such that 
$q_n\in T_{\varepsilon_n}(\Sigma)$ and $q_n\in\cC'_1\cup\cD$ for all $n$. In particular, for $n$ large enough, $q_n$ belongs to the union 
of $\cD$ and the part of $\cC_1$ between the planes $\{y_2=0\}$ and $\{y_2=1\}$, which is compact. Hence, up to extraction of a subsequence, we 
can assume that $q_n$ converges to a point $q$. We necessarily have $q\in P$, and, since $\Sigma$ is proper, $q\in\Sigma$.
 This contradicts the fact  that $\Sigma\cap P=\emptyset$, which proves the claim. 

By proposition \ref{limitinfinity}, $\cC'_{\alpha}$ converges smoothly, away from $0$, to $P\setminus\{0\}$ when 
$\alpha\to+\infty$.  Hence, for $\alpha$ large enough, $\cC'_{\alpha}\cap T_{\varepsilon}(\Sigma)\neq\emptyset$. 
Also, by continuity of the family $(\cC_{\alpha})$, we have $\cC'_{\alpha}\cap T_{\varepsilon}(\Sigma)=\emptyset$ for 
$\alpha$ close enough to $1$.

Let $\Gamma=\{\alpha\geqslant 1;\cC'_{\alpha}\cap T_{\varepsilon}(\Sigma)\neq\emptyset\}$ and $\gamma=\inf\Gamma$. 
We have $\gamma>1$. We claim that $\gamma\in\Gamma$.

If $\gamma$ is an isolated point, then it is clear. We now assume that $\gamma$ is not isolated. Then there exists a decreasing sequence 
$(\alpha_n)$ converging do $\Gamma$ and a sequence of points $(p_n)$ such that $p_n\in\cC'_{\alpha_n}\cap T_{\varepsilon}(\Sigma)$. 
We can write $p_n=X_{\alpha_n}(u_n,A_n)$ with $u_n\in[-U_{\alpha_n},U_{\alpha_n}]$ and $A_n\in\R$. 
We have $0\leqslant(y_2)_{\alpha_n}(p_n)\leqslant\varepsilon$, i.e., $$0\leqslant\frac{C_{\alpha_n}}{\alpha_n}A_n-
\frac{C_{\alpha_n}}{\alpha_n}\beta_{\alpha_n}(u_n)
-G_{\alpha_n}(u_n)\leqslant\varepsilon.$$
Since for all $n$ we have $\alpha_n\in[\gamma,\alpha_0]$, $u_n$ is bounded and so 
$|\frac{C_{\alpha_n}}{\alpha_n}\beta_{\alpha_n}(u_n)+G_{\alpha_n}(u_n)|$ is also bounded; 
moreover $C_{\alpha_n}$ is bounded from below by a positive constant. From this we deduce that $A_n$ is bounded. Consequently,
up to extraction of a subsequence, we can assume that $(u_n,A_n)$ converges to some $(u,A)\in\R$. Then, by continuity, 
$p_n$ converges to a point lying in $T_{\varepsilon}(\Sigma)$ and in $\cC'_{\gamma}$. This finishes proving the claim.

Thus there exists a point $p\in\cC'_{\gamma}\cap T_{\varepsilon}(\Sigma)$. Since $\partial\cC'_{\gamma}
\subset\cD$ (by construction of $\cD$) and $T_{\varepsilon}(\Sigma)\cap\cD=\emptyset$, $p$ is an interior point of $\cC'_{\gamma}$. 
Moreover, since $\cC'_{\alpha}\cap T_{\varepsilon}(\Sigma)=\emptyset$ for all $\alpha<\gamma$, $\cC'_{\gamma}$ lies on one 
side of $T_{\varepsilon}(\Sigma)$ in a neighbourhood of $p$. Then, by the maximum principle we get 
$T_{\varepsilon}(\Sigma)=\cC_{\gamma}$, which gives a contradiction since a horizontal catenoid is not contained in a half-space.
\end{proof}

\begin{rem}
Apart from the fact that horizontal catenoids converge to a punctured vertical plane, the key fact in this proof 
is that horizontal catenoids meet all vertical planes $\{y_2=c\}$ for $c\in\R$ (propostion \ref{convex}). This ensures that the 
sequence $(p_n)$ is bounded.

For example, for minimal surfaces in $\h\times\R$, rotational catenoids have finite height (see \cite{nelli}), and so there 
is no half-space theorems with respect to horizontal planes.
\end{rem}

\begin{rem}
Abresch and Rosenberg \cite{AR2} proved a half-space theorem with respect to surfaces of equation $x_3=c$ ($c\in\R$). It relies on the fact that rotational catenoids converge to such a surface.
\end{rem}
\section{Helicoidal minimal surfaces}\label{helicoid}
In this section we investigate the minimal surface $\cH_{\alpha}$ in $\nil$ whose Gauss map is the function $g$ defined in section \ref{harmonic},
with $\theta =0$. The calculus of Proposition \ref{integration} still holds with $C=0$, $\beta=0$ and $A=\alpha v$.

The coordinates $(y_1,y_2,y_3)$ of $\cH_{\alpha}$ are
\begin{eqnarray*}
y_1 & = & \frac{G'(u)}{\alpha}\cos\varphi(u)\sinh (\alpha v), \\
y_2 & = & -G(u), \\
y_3 & = & -\frac{G'(u)}{\alpha}\sin\varphi(u)\sinh (\alpha v).
\end{eqnarray*}
 
In particular we have $\frac{y_{3}}{y_{1}}=-\tan \vp (u)$ and $y_{2}$ only depends on $u$. This means that the intersection
of $\cH_{\alpha}$ and any vertical plane $\{y_{2}=c\}$ is a straight line. Moreover the surface is simply periodic since $\vp$ is periodic.

\section{A family of horizontal CMC $\frac12$ annuli}\label{cmc}

In this section we integrate the equations of Fernandez and Mira to construct a one-parameter
family of horizontal annuli. These surfaces are the sister
surfaces of the helicoidal surfaces construct in section \ref{helicoid}.

We consider the harmonic map $g:\Sigma \to \h$ given in section 
\ref{harmonic}, with $\theta=0$ and $g_{*}$ the conjugate harmonic map 
(see section \ref{harmonic}).

\begin{lemma} The following harmonic maps are conjugate:
$$g=\frac{\sin \vp (u) + i\sinh (\alpha v)}{\cos \vp (u) + \cosh (\alpha v)},$$
$$g_{*}=\frac{\sin \vp_{*} (u) + i\sinh (\alpha_{*} v)}{\cos \vp _{*}(u) + \cosh (\alpha_{*} v)},$$
with $\vp '^2 -\alpha^2=\cos ^2 \vp , \vp (0)=0$ and $\vp_{*}'^2-\alpha_{*} ^2=-\cos ^2 \vp_{*}, \vp_{*}(0)=0$ and $\alpha_{*}^2=\alpha ^2 +1$.  
\end{lemma}

\begin{proof}
We remark that $g_{*}$ is harmonic as $g$ in section \ref{harmonic}. Moreover $Q(g_{*})=-\frac14 \rmd z^2$.
The conformal minimal immersions in $\h \times  \R$ are given by $Y(u,v)=(g(u,v), v)$ and $Y_{*}=(g_{*}(u,v),u)$. 
To be isometric, it suffices to check that $\cosh \omega=\cosh \omega_{*}$, i.e.,
$$\cosh ^2 \omega = \frac{4|g_{u}|^2}{(1-|g|^2)^2}=1+\frac{4|g_{v}|^2}{(1-|g|^2)^2}=
\frac{4|g_{*v}|^2}{(1-|g_{*}|^2)^2}=1+\frac{4|g_{*u}|^2}{(1-|g_{*}|^2)^2}$$
These relations are equivalent to
$$\frac{\vp '^2}{\cos ^{2}\vp}=1+ \frac{\alpha ^2}{\cos ^2 \vp}=1+ \frac{\vp'^2_{*} }{\cos ^2 \vp_{*}}
=\frac{\alpha_{*}^2}{\cos ^{2}\vp _{*}}.$$
A straighforward computation shows the functions $\frac{\vp '}{\cos \vp}$ and $\frac{\alpha_{*}}{\cos \vp _{*}}$
are both solutions of
$$A'^2=(A^2-1)(A^2-\alpha^2 -1).$$
Moreover we have $\vp (0)=\vp_{*}(0)=0$ because $\alpha_{*}^2=\alpha^2+1$. This concludes the proof.
\end{proof}

In summary we have $Q(g)=\frac14\rmd z^2=-Q(g_{*})$ and $\tau =\tau^*=e^{2\omega}$.
The map $g$ induces locally a minimal graph in $\nil$ by theorem \ref{weierstrass}, with metric 
$$\lambda = \frac{\tau}{\nu^2} |\rmd z|^2= 16 \frac{ (1+|g|^2)^2}{(1-|g|^2)^4}|g_{z}|^2 |\rmd z|^2.$$
Then
\begin{eqnarray}
\label{eq:utau}
\nu^2=\frac{1-|g|^2}{1+|g|^2}=\frac{\cos ^2 \varphi}{\cosh ^2(\alpha v)}    ,\quad
 \tau=  \frac{(\varphi  ' +\alpha)^2}{\cos ^2 \varphi}=\frac{(\varphi _{*}  ' +\alpha_{*})^2}{\cos ^2 \varphi _{*}}.
\end{eqnarray}

This minimal multigraph is isometric to an immersed  CMC $\frac12$ surface 
in $\h\times \R$ with harmonic Gauss map $g_{*}:\Sigma \to \h$ 
admitting data $(-Q, \tau)$. 
For  $a_{0}\in \C$, there is a unique solution $h^*$ of the following system
\[
\left\{\begin{array}{ll}
h^*_{zz}=(\log \tau)_{z} h^*_{z} + Q\sqrt{\frac{\tau + 4 |h^*_{z}|^2}{\tau}}\\
\\
h^*_{z\bar z}=\frac{1}{4}\sqrt{\tau (\tau + 4 |h^*_{z}|^2)}\\
\\
h^*_{z}(z_{0})=a_{0}
\end{array}\right.
\]
with $\tau + 4|h^*_{z}|^2=\lambda$, and using (\ref{eq:utau}), $\vp ''+\sin \vp \cos \vp=0$, $Q=\frac14$ we obtain
\[
\left\{\begin{array}{ll}
h^*_{zz}=\alpha \tan \vp h^*_{z} + \frac{\cosh(\alpha v)}{4 \cos \vp}\\
\\
h^*_{z\bar z}=\frac{(\vp '+\alpha)^2 \cosh (\alpha v)}{4 \cos ^3 \vp}
\end{array}\right.
\]
Now set $H=h^*_{z}$, then
\[
\left\{\begin{array}{ll}
H_{z}=\alpha \tan \vp H + \frac{\cosh(\alpha v)}{4\cos \vp}\\
\\
H_{\bar z}=\frac{(\vp '+\alpha)^2 \cosh (\alpha v)}{4 \cos ^3 \vp}\\
\\
H (z_{0})=a_{0}.
\end{array}\right.
\]
Then
$$H(u,v)=\frac{\cos \vp}{2(\alpha-\vp ')}\left(i\sinh (\alpha v)-\tan \vp \cosh (\alpha v)     \right)+K_{1}(u) e^{i (\alpha \tan \vp)v}  
+ K_{2}(u)$$
where $K_{i}'=a \tan \vp K_{i}$ and $K_{1}(u_{0}),K_{2}(u_{0})$ are chosen to have $H(z_{0})=a_{0}$. 
It is a two-parameter family and we are interested in $a_{0}$ such that, $K_{1}=K_{2}=0$, i.e.,
the solution with $\tau +4|H|^2=\lambda$. The solution is periodic in $u$ and by \eqref{eq:utau} we have
$$h^*=\frac{ \cos \vp \cosh (\alpha v)}{\alpha  (\vp'-\alpha)}=\frac{\cos \vp_{*} \cosh (\alpha v)}{\alpha (\vp_{*}'-\alpha_{*})}.$$
Now we consider $$(G_{1},G_{2},G_{3})=\left( \frac{2g_{*}}{1-|g_{*}|^2}, \frac{1+|g_{*}|^2}{1-|g_{*}|^2}\right)=
\left(\tan \vp_{*},\frac{\sinh (\alpha_{*} v)}{\cos \vp_{*}},\frac{\cosh (\alpha_{*} v)}{\cos \vp_{*}}\right)$$
and we compute the horizontal component $F^*=\frac{X_{1}+iX_{2}}{1+X_{3}}$ given by
$$X_{j}=\frac{8{\rm Re}\left(G_{j,z}(4\bar Q h_{z} +\tau h_{\bar z}) \right)}{\tau ^2 -16 |Q|^2}+
G_{j}\sqrt{\frac{\tau + 4 |h_{z}|^2}{\tau}}$$
i.e.,
$$X_{j}=(\alpha_{*} - \vp'_{*})\left(\frac{G_{j,u}h_{u}}{\vp '_{*}}+\frac{G_{j,v}h_{v}}{\alpha_{*}}\right)+
G_{j}\frac{\cosh (\alpha v)}{\cos \vp}.$$
Straightforward computations give
\[
\begin{array}{ll}
h_{u}=\frac{\alpha_{*} \sin \vp_{*} \cosh (\alpha v)}{\alpha (\vp '_{*}-\alpha_{*})} &
 h_{v}=\frac{\cos  \vp_{*}\sinh (\alpha v)}{(\vp'_{*}-\alpha_{*})} \\
&\\
G_{1,u}=\frac{\vp '_{*}}{\cos ^2 \vp_{*}} & G_{1,v}=0\\
& \\
G_{2,u}=\frac{\vp'_{*}\sin \vp_{*}\sinh (\alpha_{*} v)}{\cos ^2 \vp_{*}} &
G_{2,v}=\frac{\alpha_{*}\cosh (\alpha_{*} v)}{\cos  \vp_{*}}\\
& \\
G_{3,u}=\frac{\vp'_{*}\sin \vp_{*}\cosh (\alpha_{*} v)}{\cos ^2 \vp_{*}} &
G_{3,v}=\frac{\alpha_{*}\sinh (\alpha_{*} v)}{\cos  \vp_{*}}
\end{array}
\]
Inserting the explicit value above, setting $f(u)=
\frac{\alpha \cos \vp_{*}-\alpha _{*} \cos \vp}{\alpha \cos \vp \cos ^2 \vp_{*}}$, we obtain
\[
\left\{\begin{array}{l}
\displaystyle{X_{1}=\cosh (\alpha v)\sin \vp_{*}(u) f(u)} \\
\\
\displaystyle{X_{2}=\cosh (\alpha v) \sinh (\alpha_{*} v)\left (f(u) +\frac{\alpha _{*}}{\alpha} \right)-
\cosh ( \alpha_{*} v) \sinh (\alpha v) } \\
\\
\displaystyle{X_{3}=\cosh (\alpha v) \cosh (\alpha_{*} v) \left( f(u)+\frac{\alpha _{*}}{\alpha} \right)-
\sinh ( \alpha_{*} v) \sinh (\alpha v)}
\end{array}\right.
\]

Now we are interested in level curves at height zero. Then by \eqref{U2}, we have
 $\vp (- U/2)=\vp_{*}(-U/2)=+\pi/2$ and
$\vp (U/2)=\vp_{*}(U/2)=-\pi/2$. We have $h(-U/2,v)=h(U/2,v)=0$. Using the ODE of
$\vp$ and $\vp_{*}$ we obtain $f(u) \to \gamma= \frac{-1}{2\alpha  \alpha _{*}}$
when $u \to \pm U/2$. Then the horizontal curve $h^*=0$ has two connected components given 
by $F^*(\pm U/2,v)=\frac{X_{1}+iX_{2}}{1+X_{3}}$ with
\[
\left\{\begin{array}{l}
\displaystyle{X_{1}=\cosh (\alpha v)\sin \vp_{*}(\pm U/2) \gamma }\\
\\
\displaystyle{X_{2}=\cosh (\alpha v) \sinh (\alpha_{*} v)\left (\gamma +\frac{\alpha_*}{\alpha} \right)-
\cosh ( \alpha_{*} v) \sinh (\alpha v) } \\
\\
\displaystyle{X_{3}=\cosh (\alpha v) \cosh (\alpha_{*} v) \left( \gamma+\frac{\alpha_*}{\alpha} \right)-
\sinh ( \alpha_{*} v) \sinh (\alpha v)}
\end{array}\right.
\]
We remark that $\gamma +\frac{\alpha_{*}}{\alpha}=\frac{2\alpha ^2 +1}{2\alpha \alpha_*} \geqslant 1$. 
 In the disk model $F^* (- U/2,v)$ and $F^*(U/2,v)$ are symetric with respect to the $y$-axis of the disk
 ($X_{1}(-U/2,v)=-X_{1}(U/2,v)$). Then we will study $F^*(-U/2,v)$. It is a curve linking
 the point $(0,-1)$ to $(0,1)$ in the unit disk and staying in $\re F^* >0$. We prove that this curve is embedded
 and behaves like a generatrix of a Bryant catenoid in hyperbolic three-space.
In the half-plane model of $\h$, this curve is given by 
\begin{eqnarray*}
\tilde F (-U/2,v) & = & \left( \frac{X_{1}}{X_{3}-X_{2}},\frac{1}{X_{3}-X_{2}}\right) \\
& = & \left(\frac{\gamma e^{\alpha_{*}v}}{(\gamma + \frac{\alpha_*}{\alpha} +\tanh(\alpha v) ) },
 \frac{e^{\alpha_{*}v}}{\cosh (\alpha v)(\gamma + \frac{\alpha_*}{\alpha} +\tanh(\alpha v)) } \right).
\end{eqnarray*}
By a straighforward computation we can see that the map
 $$v \mapsto \frac{\gamma e^{\alpha_{*}v}}{(\gamma + \frac{\alpha_*}{\alpha} +\tanh(\alpha v)) } $$
is strictly increasing for $\alpha \leqslant 1$ and has exactly one point where the derivative is zero when $\alpha =1$.
 When $\alpha \geqslant 1$ the function has two local extrema at points $v$ where 
 $\tanh (\alpha v_{\pm})=\frac{-\sqrt{\alpha ^2 +1}\pm \sqrt{\alpha ^2-1}}{2 \alpha}) <0$. 
 For $v >0$,  $v(x_{1})$ is a well defined function and 
 $\frac{1}{X_{3}-X_{2}}=x_{1}^{1-\frac{\alpha}{\alpha_{*}}}q(x_{1})$ 
 where $q(x_{1})$ is a bounded function having a positive limit at infinity. 
 When $\alpha \to \infty$, the curve converges to two tangent horocycle.
 
 The immersion for $u\in [-U/2,U/2]$ is a graph over a simply connected domain of $\h$. We complete
 it by reflection about the horizontal plane of height zero in $\h \times \R$ to obtain a properly embedded annulus.

 
\section{Appendix: proofs of lemmas \ref{lemmaL} and lemma \ref{limitthetainfinity}} \label{appendix}

We use the notations of section \ref{harmonic}. We notice that, for $(\alpha,\theta)\in\Omega$,  $\frac{-\alpha\sqrt{P_{\alpha,\theta}(x)}+\alpha^2}{x^2}$ 
can be extended smoothly at $x=0$ and that
$$L(\alpha,\theta)=\int_{-1}^1\frac{-\alpha\sqrt{P_{\alpha,\theta}(x)}+\alpha^2+C_{\alpha,\theta}^2x^4}
{x^2\sqrt{(1-x^2)}\sqrt{P_{\alpha,\theta}(x)}}.$$

\begin{lemmaintro}[lemma \ref{lemmaL}]
Let $\alpha>0$. Then there exists a unique $\tilde\theta_{\alpha}\in(0,\theta^+_{\alpha})\cap(0,\frac\pi4)$ such that 
$$L(\alpha,\tilde\theta_{\alpha})=0.$$
\end{lemmaintro}

\begin{proof}
We have $$L(\alpha,\theta)=\int_{-1}^1\frac{l(\alpha,\theta,x)}{\sqrt{1-x^2}}\rmd x
=\int_{-\frac{\pi}2}^{\frac{\pi}2}l(\alpha,\theta,\sin t)\rmd t$$ where $l$ is a smooth function on $\Omega\times[-1,1]$. Hence $L$ is smooth on $\Omega$.

We have $L(\alpha,0)<0$, since $C_{\alpha,0}=0$.

We first deal with the case where $\alpha>\frac1{\sqrt{2}}$, i.e., $\theta^+_{\alpha}>\frac\pi4$. Since the integrand in $L(\alpha,\frac\pi4)$ is positive for all $x\in(-1,1)$, we have $L(\alpha,\frac\pi4)>0$. Hence, by continuity, there exists $\tilde\theta_{\alpha}\in(0,\frac\pi4)$ such that $L(\alpha,\tilde\theta_{\alpha})=0$.

We now deal with the case where $\alpha\leqslant\frac1{\sqrt{2}}$, i.e., $\theta^+_{\alpha}\leqslant\frac\pi4$. We have $C_{\alpha,\theta}^2\to1-\alpha^2$, $\rho^-_{\alpha,\theta}\to 1$ and $\rho^+_{\alpha,\theta}\to\frac{\alpha^2}{1-\alpha^2}>0$ when $\theta\to\theta^+_{\alpha}$. We have
$$L(\alpha,\theta)=L_1(\alpha,\theta)+L_2(\alpha,\theta)$$ with
$$L_1(\alpha,\theta)=\int_{-1}^1-\frac{2\alpha C_{\alpha,\theta}^2}
{\alpha+\sqrt{P_{\alpha,\theta}(x)}}\sqrt{\frac{1-x^2}{P_{\alpha,\theta}(x)}}\rmd x,$$
$$L_2(\alpha,\theta)=\int_{-1}^1\frac{2\alpha C_{\alpha,\theta}^2-\alpha\cos(2\theta)+C_{\alpha,\theta}^2x^2\sqrt{P_{\alpha,\theta}(x)}}
{\sqrt{(1-x^2)P_{\alpha,\theta}(x)}(\alpha+\sqrt{P_{\alpha,\theta}(x)})}\rmd x.$$
We claim that $\frac{1-x^2}{P_{\alpha,\theta}(x)}$ is uniformly bounded (in $x$) when $\theta\to\theta^+_{\alpha}$. Indeed we have 
\begin{eqnarray*}
\left|\frac{1-x^2}{P_{\alpha,\theta}(x)}-\frac1{\alpha^2+(1-\alpha^2)x^2}\right| & = &
\left|\frac{(1-2\alpha^2-\cos(2\theta))x^2-(1-\alpha^2-C_{\alpha,\theta}^2)x^4}
{C_{\alpha,\theta}^2(\rho^-_{\alpha,\theta}-x^2)(\rho^+_{\alpha,\theta}+x^2)(\alpha^2+(1-\alpha^2)x^2)}\right|
\\
& \leqslant & 
\frac{|1-2\alpha^2-\cos(2\theta)|+|1-\alpha^2-C_{\alpha,\theta}^2|}
{C_{\alpha,\theta}^2(\rho^-_{\alpha,\theta}-1)\rho^+_{\alpha,\theta}\alpha^2} \\
& \leqslant &
\frac{(\cos(2\theta)+1+2\alpha^2)(1-\cos(2\theta))}{4\alpha^2C_{\alpha,\theta}^2\rho^+_{\alpha,\theta}}.
\end{eqnarray*}
This upper bound has a finite limit when $\theta\to\theta^+_{\alpha}$. This proves the claim. Consequently, $L_1(\alpha,\theta)$ is bounded when $\theta\to\theta^+_{\alpha}$.
Moreover we have $2\alpha C_{\alpha,\theta}^2-\alpha\cos(2\theta)\to\alpha$ when $\theta\to\theta^+_{\alpha}$, so there exists a positive constant $c_{\alpha}$ such that, for $\theta$ close enough to $\theta^+_{\alpha}$,
$$L_2(\alpha,\theta)\geqslant\int_{-1}^1\frac{c_{\alpha}}{\sqrt{(1-x^2)(\rho^-_{\alpha,\theta}-x^2)}}.$$
Since $\rho^-_{\alpha,\theta}\to 1$ when $\theta\to\theta^+_{\alpha}$, we obtain that $L_2(\alpha,\theta)\to+\infty$ when $\theta\to\theta^+_{\alpha}$. We conclude that $L(\alpha,\theta)\to+\infty$ when $\theta\to\theta^+_{\alpha}$. Hence, by continuity, there exists $\tilde\theta_{\alpha}\in(0,\theta_0)$ such that $L(\alpha,\tilde\theta_{\alpha})=0$.

To prove the uniqueness of $\tilde\theta_{\alpha}$, it suffices to prove that $\theta\mapsto L(\alpha,\theta)$ is increasing on $(0,\theta^+_{\alpha})\cap(0,\frac\pi4)$. A straightforward computation gives
$$\frac{\partial}{\partial\theta}\left(\frac{-\alpha\sqrt{P_{\alpha,\theta}(x)}+\alpha^2+C_{\alpha,\theta}^2x^4}
{x^2\sqrt{P_{\alpha,\theta}(x)}}\right)
=\frac{C_{\alpha,\theta}}{\alpha}\frac{B_1}{P_{\alpha,\theta}(x)^{\frac32}}$$
with $$B_1=2\cos(2\theta)x^2P_{\alpha,\theta}(x)
+(2\alpha^2+\cos(2\theta)x^2)(\alpha^2+C_{\alpha,\theta}^2x^4)>0.$$



This shows that $L(\alpha,\theta)$ is increasing in $\theta$, which proves the uniqueness of $\tilde\theta_{\alpha}$.
\end{proof}

\begin{lemmaintro}[lemma \ref{limitthetainfinity}]
We have $$\lim_{\alpha\to+\infty}\tilde\theta_{\alpha}=\frac{\pi}4.$$
\end{lemmaintro}

\begin{proof}
We set $C_{\alpha}=C_{\alpha,\tilde\theta_{\alpha}}$ and $P_{\alpha}(x)=P_{\alpha,\tilde\theta_{\alpha}}(x)$.
We first notice that $C_{\alpha}\leqslant\frac1{2\alpha}$ and, for $\alpha\geqslant\frac12$ and $x\in[-1,1]$, $\alpha^2-1\leqslant P_{\alpha}(x)\leqslant\alpha^2+1$.

We have $$0=\alpha L(\alpha,\tilde\theta_{\alpha})=I_1(\alpha)-\cos(2\tilde\theta_{\alpha})I_2(\alpha)
+I_3(\alpha)$$
with
$$I_1(\alpha)=\int_{-1}^1\frac{2\alpha^2C_{\alpha}^2}{\sqrt{P_{\alpha}(x)}(\alpha+\sqrt{P_{\alpha}(x)})}
\frac{x^2\rmd x}{\sqrt{1-x^2}}=\rmO\left(\frac1{\alpha^2}\right),$$
$$I_2(\alpha)=\int_{-1}^1\frac{\alpha^2}{\sqrt{P_{\alpha}(x)}(\alpha+\sqrt{P_{\alpha}(x)})}
\frac{\rmd x}{\sqrt{1-x^2}}\geqslant\frac{\pi\alpha^2}{\sqrt{\alpha^2+1}(\alpha+\sqrt{\alpha^2+1})},$$
$$I_3(\alpha)=\int_{-1}^1\frac{\alpha C_{\alpha}^2}{\alpha+\sqrt{P_{\alpha}(x)}}
\frac{x^2\rmd x}{\sqrt{1-x^2}}=\rmO\left(\frac1{\alpha^2}\right).$$
Hence $\cos(2\tilde\theta_{\alpha})\to 0$, which proves the lemma.
\end{proof}

\bibliographystyle{alpha}
\bibliography{catenoid12}

\begin{thebibliography}{HSET05}

\bibitem[AR04]{AR1}
U.~Abresch and H.~Rosenberg.
\newblock A {H}opf differential for constant mean curvature surfaces in {${\bf
  S}\sp 2\times{\bf R}$} and {${\bf H}\sp 2\times{\bf R}$}.
\newblock {\em Acta Math.}, 193(2):141--174, 2004.

\bibitem[AR05]{AR2}
U.~Abresch and H.~Rosenberg.
\newblock Generalized {H}opf differentials.
\newblock {\em Mat. Contemp.}, 28:1--28, 2005.

\bibitem[Bry87]{bryant}
R.~Bryant.
\newblock Surfaces of mean curvature one in hyperbolic space.
\newblock {\em Ast\'erisque}, (154-155):12, 321--347, 353 (1988), 1987.
\newblock Th\'eorie des vari\'et\'es minimales et applications (Palaiseau,
  1983--1984).

\bibitem[CM99]{col}
T.~H. Colding and W.~P. Minicozzi, II.
\newblock {\em Minimal surfaces}, volume~4 of {\em Courant Lecture Notes in
  Mathematics}.
\newblock New York University Courant Institute of Mathematical Sciences, New
  York, 1999.

\bibitem[CR07]{cr}
P.~Collin and H.~Rosenberg.
\newblock The {J}enkins-{S}errin theorem for minimal graphs in homogeneous
  $3$-manifolds.
\newblock Preprint in preparation, 2007.

\bibitem[CY75]{chengyau2}
S.~Y. Cheng and S.~T. Yau.
\newblock Differential equations on {R}iemannian manifolds and their geometric
  applications.
\newblock {\em Comm. Pure Appl. Math.}, 28(3):333--354, 1975.

\bibitem[Dan04]{daniel}
B.~Daniel.
\newblock Isometric immersions into $\mathbb{S}^n \times \mathbb{R}$ and
  $\mathbb{H}^n \times \mathbb{R}$ and applications to minimal surfaces.
\newblock arXiv:math.DG/0406426, to appear in {\it Trans. Amer. Math. Soc.},
  2004.

\bibitem[Dan06]{dan2}
B.~Daniel.
\newblock The {G}auss map of minimal surfaces in the {H}eisenberg group.
\newblock arXiv:math.DG/0606299, 2006.

\bibitem[Dan07]{dan}
B.~Daniel.
\newblock Isometric immersions into 3-dimensional homogeneous manifolds.
\newblock {\em Comment. Math. Helv.}, 82(1):87--131, 2007.

\bibitem[FM07a]{fermira1}
I.~Fern{\'a}ndez and P.~Mira.
\newblock Harmonic maps and constant mean curvature surfaces in {$\mathbb{H}\sp
  2\times\mathbb{R}$}.
\newblock {\em Amer. J. Math.}, 129(4):1145--1181, 2007.

\bibitem[FM07b]{fermira2}
I.~Fern\'andez and P.~Mira.
\newblock Holomorphic quadratic differentials and the {B}ernstein problem in
  {H}eisenberg space.
\newblock arXiv:0705.1436, to appear in {\it Trans. Amer. Math. Soc.}, 2007.

\bibitem[FMP99]{mercuri}
C.~Figueroa, F.~Mercuri, and R.~Pedrosa.
\newblock Invariant surfaces of the {H}eisenberg groups.
\newblock {\em Ann. Mat. Pura Appl. (4)}, 177:173--194, 1999.

\bibitem[Hau06]{hau}
L.~Hauswirth.
\newblock Minimal surfaces of {R}iemann type in three-dimensional product
  manifolds.
\newblock {\em Pacific J. Math.}, 224(1):91--117, 2006.

\bibitem[HM90]{HM}
D.~Hoffman and W.~H. Meeks, III.
\newblock The strong halfspace theorem for minimal surfaces.
\newblock {\em Invent. Math.}, 101(2):373--377, 1990.

\bibitem[HR07]{rosenhau}
L.~Hauswirth and H.~Rosenberg.
\newblock Minimal surfaces of finite total curvature in $\mathbb{H} \times
  \mathbb{R}$.
\newblock Preprint, http://perso-math.univ-mlv.fr/users/hauswirth.laurent, to
  appear in {\it Mat. Contemp.}, 2007.

\bibitem[HRS07]{hrs2}
L.~Hauswirth, H.~Rosenberg, and J.~Spruck.
\newblock On complete mean curvature $\frac12$ surfaces in
  $\mathbb{H}^2\times\mathbb{R}$.
\newblock Preprint, http://perso-math.univ-mlv.fr/users/hauswirth.laurent,
  2007.

\bibitem[HSET05]{touhau}
L.~Hauswirth, R.~S\'a~Earp, and E.~Toubiana.
\newblock Associate and conjugate minimal immersion in $\mathbb{M}\times
  \mathbb{R}$.
\newblock arXiv:math.DG/0512112, to appear in {\it Tohoku Math. J.}, 2005.

\bibitem[NR02]{nelli}
B.~Nelli and H.~Rosenberg.
\newblock Minimal surfaces in {${\mathbb H}\sp 2\times{\mathbb R}$}.
\newblock {\em Bull. Braz. Math. Soc. (N.S.)}, 33(2):263--292, 2002.

\bibitem[UY93]{uy}
M.~Umehara and K.~Yamada.
\newblock Complete surfaces of constant mean curvature {$1$} in the hyperbolic
  {$3$}-space.
\newblock {\em Ann. of Math. (2)}, 137(3):611--638, 1993.

\bibitem[WA94]{wan-au}
T.~Y.-H. Wan and T.~K.-K. Au.
\newblock Parabolic constant mean curvature spacelike surfaces.
\newblock {\em Proc. Amer. Math. Soc.}, 120(2):559--564, 1994.

\bibitem[Wan92]{wan}
T.~Y.-H. Wan.
\newblock Constant mean curvature surface, harmonic maps, and universal
  {T}eichm\"uller space.
\newblock {\em J. Differential Geom.}, 35(3):643--657, 1992.

\end{thebibliography}

\end{document}